\documentclass[11pt,a4paper]{article}
\usepackage{amsmath,amssymb,amsthm,xspace,a4wide}
\usepackage{amscd,latexsym,amscd,stmaryrd,xy}
\xyoption{all}
\usepackage{enumerate}
\numberwithin{equation}{section}

\newtheorem{satz}{Theorem}%C
\newtheorem{lemma}{Lemma}[section]
\newtheorem{remark}[lemma]{Remark}%[section]
\newtheorem{theorem}[satz]{Theorem}%[section]
\newtheorem{question}[lemma]{Question}%[section]

\newtheorem{claim}[lemma]{Claim}%[section]
\newtheorem{proposition}[lemma]{Proposition}%[section]
\newcommand{\PT}{I^\Delta}

\newcommand{\Z}{\ensuremath{\mathbb{Z}}}

\newcommand{\Oo}{\ensuremath{\mathcal{O}}}

\newcommand{\g}{\ensuremath{\mathfrak{g}}}

\newcommand{\h}{\ensuremath{\mathfrak{h}}}
\newcommand{\n}{\ensuremath{\mathfrak{n}}}

\newcommand{\tto}{\twoheadrightarrow}

\renewcommand{\phi}{\varphi}

\newcommand{\Hom}{\operatorname{Hom}}
\newcommand{\End}{\operatorname{End}}

\newcommand{\Ker}{\operatorname{Ker}}

\renewcommand{\Im}{\operatorname{Im}}
\newcommand{\Ext}{\operatorname{Ext}}

\newcommand{\ho}{\mathtt{H}}
\newcommand{\erem}{\hfill$\blacksquare$}

\newcommand{\cA}{{\mathcal A}}
\newcommand{\cB}{{\mathcal B}}
\newcommand{\cC}{{\mathcal C}}

\newcommand{\cH}{{\mathcal H}}

\newcommand{\cL}{{\mathcal L}}

\newcommand{\cO}{{\mathcal O}}
\newcommand{\cP}{{\mathcal P}}
\newcommand{\cR}{{\mathcal R}}
\newcommand{\cS}{{\mathcal S}}
\newcommand{\cT}{{\mathcal T}}

\newcommand{\cZ}{{\mathcal Z}}

\newcommand{\Coinv}{\cC}

\newcommand{\mg}{\mathfrak{g}}
\newcommand{\mh}{\mathfrak{h}}

\newcommand{\mV}{\mathbb{V}}

\newcommand{\mC}{\mathbb{C}}

\newcommand{\mZ}{\mathbb{Z}}
\newcommand{\SL}{\mathfrak{sl}}

\newcommand{\la}{\lambda}

\newcommand{\HOM}{\operatorname{Hom}}

\newcommand{\EXT}{\operatorname{Ext}}

\newcommand{\END}{\operatorname{End}}

\newcommand{\MOF}{\operatorname{mod}}

\newcommand{\op}{\operatorname}

\newcommand{\KER}{\operatorname{ker}}
\newcommand{\ADJ}{\operatorname{adj}}

\newcommand{\IM}{\operatorname{im}}

\newcommand{\inj}{\hookrightarrow}
\newcommand{\surj}{\mbox{$\rightarrow\!\!\!\!\!\rightarrow$}}

\newcommand{\XX}{\mathrm{X}}
\newcommand{\YY}{\mathrm{Y}}
\newcommand{\ZZ}{\mathrm{Z}}
\newcommand{\CC}{\mathrm{C}}
\newcommand{\KK}{\mathrm{K}}
\newcommand{\TT}{\mathrm{T}}
\newcommand{\GG}{\mathrm{G}}
\newcommand{\FF}{\mathrm{Q}}

\newcommand{\wZZ}{\hat{\mathrm{Z}}}
\newcommand{\EE}{\mathrm{E}}
\newcommand{\ID}{\mathrm{ID}}
\newcommand{\rd}{\mathrm{d}}
\newcommand{\rg}{\mathrm{g}}
\newcommand{\rt}{\mathrm{g}'}
\newcommand{\ri}{\mathrm{i}}

\newcommand{\rp}{\mathrm{p}}

\newcommand{\rz}{\mathrm{z}}
\newcommand{\rk}{\mathrm{m}}

\newcommand{\raa}{\mathrm{a}}

\newcommand{\rh}{\mathrm{h}}

\newcommand{\rc}{\mathrm{c}}
\newcommand{\rkw}{\mathrm{k}}

\begin{document}
\title{On functors associated to a simple root}
\author{Volodymyr Mazorchuk and Catharina Stroppel}
\date{}

\maketitle
\begin{abstract}
Associated to a simple root of a finite-dimensional complex semisimple
 Lie algebra, there are several endofunctors
(defined by Arkhipov, Enright, Frenkel, Irving, Jantzen, Joseph, Mathieu,
 Vogan and Zuckerman) on the BGG category $\Oo$. We study their relations,
compute cohomologies of their derived functors and describe the
monoid generated by Arkhipov's and  Joseph's functors and the monoid
generated by Irving's functors. It turns out that the
endomorphism rings of all elements in these monoids are isomorphic.
We prove that the functors give rise to an action of the
 singular braid monoid on the bounded derived category of $\Oo_0$.
 We also use Arkhipov's, Joseph's and Irving's functors to produce
 new  generalized tilting modules.
\end{abstract}

\section*{Introduction}\label{s0}

Braid group actions via auto-equivalences of derived categories play
an important role in different areas of mathematics and physics.
They arise naturally, for example, in algebra, algebraic geometry,
representation theory, string theory, symplectic geometry, topology etc.
In this paper we focus on endofunctors of the BGG category $\cO$,
associated to a semi-simple complex (finite-dimensional) Lie algebra,
which give rise to such braid group actions. These endofunctors
turned out to be very useful and motivating in different areas of
mathematics.  As examples one could mention for instance
\begin{itemize}
\item the Kazhdan-Lusztig combinatorics of the category $\cO$
(see e.g. \cite{KL});
\item the Serre functor for the bounded derived category of the
category of perverse sheaves on the flag variety
(see e.g. \cite{BBM});
\item structure and characters of tilting modules
(see e.g. \cite{So4});
\item finding functorial invariants of tangles and links
(see e.g. \cite{BFK});
\item defining $\mathfrak{sl}_2$-categorifications and Brou{\'e}'s conjecture
(see e.g. \cite{CR});
\item semi-infinite cohomology and Wakimoto modules
(see e.g. \cite{Ark,Ara}).
\end{itemize}
Of course, there are many more. The main properties of the functors used
in the above examples are the following: they provide a connection between
projective and tilting modules, and they define a categorification
of braid groups and Hecke algebras. To define a categorification, one
usually has to check relations between certain functors and even
between natural transformations. In practice, this is quite often
a rather cumbersome technical work (see e.g. \cite{CR}).

The aim of the present paper is the following:
\begin{itemize}
\item to collect from the wide-spread literature all the endofunctors
of the category $\cO$, associated to a simple root of the Lie algebra,
and give an insight into the interplay of these functors
(Theorem~\ref{tmain2}, Theorem~\ref{tmain1});
\item describe the non-trivial relations between these endofunctors
(Theorem~\ref{tmain3}, Theorem~\ref{tmain3new});
\item describe natural transformations between these endofunctors
(Theorem~\ref{tmain4}, Theorem~\ref{tmain4.101});
\item describe the connection between the injective, projective and
tilting modules, and construct new generalized tilting modules
(Theorem~\ref{tmain5}, Theorem~\ref{tmain6});
\item show that the structure of the category $\cO$ is not
completely determined by the Kazhdan-Lusztig combinatorics
(Remark~\ref{different}, Proposition~\ref{p2.5});
\item construct a new (unexpected) categorification of
the Baez-Birman singular braid monoid via endofunctors of the
principal block of $\cO$, which might lead to a categorification of
Vassiliev invariants (Theorem~\ref{singularbraid}).
\end{itemize}

We formulate all the results in the next section and try to hide the
technicalities (as far as possible) in the proofs, which follow afterwards.

\section{The results}\label{s1}

\subsection{Notation and the setup}\label{s1.1}

Let $\g$ be a semisimple complex finite-dimensional Lie algebra with
a fixed triangular decomposition $\g=\n_-\oplus\h\oplus\n_+$. Let $W$
be the corresponding Weyl group with the length function $l$, the unit
element $e$, the longest element $w_0$, and the Bruhat ordering $<$.
The letter $\rho$ denotes the half-sum of all positive roots.
There is the so-called dot-action of $W$
on $\mh^\ast$ defined as $w\cdot\la=w(\la+\rho)-\rho$. Let $\Oo$
denote the BGG-category $\Oo$ introduced in \cite{BGG} and $\cO_0$
its principal block, that is the indecomposable block of $\Oo$
containing the trivial $\g$-module. For a simple reflection $s$ let
$\g^s$ denote the corresponding minimal parabolic subalgebra of $\g$,
strictly containing $\mh\oplus\n_+$. We denote by $\Oo_0^s$ the
corresponding parabolic subcategory, which consists of all locally
$\g^s$-finite objects from $\Oo_0$. We call a module {$s$-free}, if
none of the composition factors in its socle is $\g^s$-finite. Let
$\Coinv=S(\mh)/(S(\mh)_+^{W\cdot})$  be the coinvariant algebra of
$W$ with respect to the dot-action. Its subalgebra of $s$-invariants
(under the usual action) is denoted by $\Coinv^s$ (see \cite[2.4]{S}).
For $x\in W$ we denote by  $\Delta(x)\in\cO_0$ the Verma module of
the highest weight $x\cdot0$ and by $P(x)$ its projective cover with
simple head $L(x)$. Associated to a fixed simple reflection $s$ we
have the following endofunctors of $\cO_0$:
\begin{itemize}
\item The {\em translation functor} $\theta=\theta_s$ {\em through the
$s$-wall} (see e.g. \cite[Section~3]{J}). Let $\Oo_s$ be a singular
integral block of $\Oo$ with stabilizer $\{e,s\}$. The functor $\theta_s$
is defined as the composition of the translation functors
$\theta_s^{on}:\cO_0\to \Oo_s$ and $\theta_s^{out}:\cO_0\to \Oo_s$
(see \cite[4.12]{Ja}), which are both left and right adjoint to each other.
In particular, the functor $\theta$ is exact and self-adjoint. For each
$x\in W$ the module $\theta_s\Delta(x)$ is a (unique up to isomorphism)
indecomposable extension of $\Delta(x)$ and $\Delta(xs)$.
\item The {\em shuffling functor} $\CC=\CC_s$ (see \cite[Section 3]{I}).
This functor is defined as follows: we fix an adjunction morphism
$\ADJ_s:\ID\rightarrow\theta$ and define $\CC$ as the cokernel of $\ADJ_s$.
Up to isomorphism the definition does not depend on the choice of $\ADJ_s$.
As both $\ID$ and $\theta$ are exact functors, the functor $\CC$ is right
exact by the Snake Lemma.
\item Dually, we have the {\em coshuffling functor} $\KK=\KK_s$. This
functor is defined as follows: we fix an adjunction morphism
$\ADJ^s:\theta\rightarrow\ID$ and define $\KK$ as the kernel of $\ADJ^s$.
Up to isomorphism the definition does not depend on the choice of $\ADJ^s$.
As both $\ID$ and $\theta$ are exact functors, the functor $\KK$ is left
exact by the Snake Lemma.
\item {\em Zuckerman's functor} $\ZZ=\ZZ_s$ (see \cite{Z}). This functor
is given by taking the maximal quotient which belongs to the category
$\Oo_0^s$. The functor $\ZZ$ is the left adjoint to the natural inclusion
of the category $\Oo_0^s$ into $\Oo_0$, in particular, $\ZZ$ is right exact.
\item {\em Joseph's completion} $\GG=\GG_s$ is defined in \cite[Section 2]{J}
in the following way: $\GG({}_-)=\mathcal{L}(\Delta(s),{}_-)\otimes \Delta(e)$,
where for $M\in \Oo_0$ the set $\mathcal{L}(\Delta(s),M)$ is a $\g$-bimodule,
consisting of all $\mathbb{C}$-linear maps from $\Delta(s)$ to $M$, which are
locally finite with respect to the adjoint action of $\g$. The functor $\GG$
is left exact.
\item {\em Arkhipov's twisting functor} $\TT=\TT_s$ is defined as the
composition of two functors (see e.g.  \cite{AS}). To define these
functors we first take a non-zero element $X_{-\alpha}$ from the negative
root subspace of $\mathfrak{g}$ associated with $s$, and localize the
universal enveloping algebra $U(\mathfrak{g})$ with respect to the powers of
$X_{-\alpha}$. We obtain the localized algebra $U_s$. The first functor is
then tensoring with the $U(\mathfrak{g})$-bimodule
$U_s/U(\mathfrak{g})$, which, in fact, does not preserve the category $\Oo_0$.
To obtain an endofunctor of $\Oo_0$ we compose
$U_s/U(\mathfrak{g})\otimes_{U(\mathfrak{g})}{}_-$ with the functor
of twisting by an automorphism of $\mathfrak{g}$, corresponding to $s$.
It turns out that the resulting functor is up to isomorphism of functors
independent of the above choice of an automorphism of $\mathfrak{g}$.
The functor $\TT$ is right exact. In fact, in
\cite[Corollary~4]{KM} it is shown that $\TT$ is left adjoint to $\GG$.
\item The functor $\FF$ given as the cokernel of a natural
transformation $\rg:\ID\rightarrow G$ as considered in \cite[2.4]{J}
(which is unique up to a non-zero scalar). The functor $\FF$ is neither
left nor right exact.
\item {\it Enright's completion functor} $\EE=\EE_s$. This functor was
originally defined in  \cite{En} via a complicated procedure, based on
the theory of highest weight $\mathfrak{sl}_2$-modules, applied to the
$\mathfrak{sl}_2$-subalgebra of $\g$ associated with $s$. The original
definition was generalized and extended in \cite{De} and \cite{J}.
In \cite[Section~4]{KM} it is shown that $\EE=\GG^2$.
\end{itemize}

The functor $\ZZ$ can be characterized as the functor of taking the
maximal quotient which is annihilated by $\TT$ (or, equivalently,
by $\GG$). We  define $\wZZ:\Oo_0\to\Oo_0$ as the endofunctor given
by taking the maximal quotient annihilated by $\CC$ (or, equivalently,
by $\KK$), i.e. the maximal quotient containing only composition
factors of the form $L(y)$, $y<ys$. Although the definition is very
similar, the properties of the functors $\ZZ$ and $\wZZ$ are quite
different (see Remark~\ref{different} and Theorem~\ref{tmain3}
below).

Let $\rd$ be the usual contravariant duality on $\Oo_0$. For an
endofunctor $\XX$ of $\Oo_0$ we denote by $\XX'$ the composition
$\XX'=\rd \XX\rd$. If $\XX_1$, $\XX_2$, $\YY$ are endofunctors on
$\Oo_0$ and $h\in \Hom(\XX_1,\XX_2)$ we denote by $h_{\YY}$ the
induced natural transformation in $\Hom(\XX_1\YY,\XX_2\YY)$.
For $h\in \Hom(\XX_1,\XX_2)$ we also set $h'=\rd\, h_{\rd}\in
\Hom(\XX_1',\XX_2')$.

In Section~\ref{s2} we give a more elegant proof of the fact
$\GG\cong \TT'$ from \cite[Theorem 4]{KM}. This result allows as to simplify
the exposition and redefine Arkhipov's functor as $\TT=\GG'$.  In
Section~\ref{s2} we also prove some similarities between the
pairs $(\TT,\GG)$ and $(\CC,\KK)$ of functors
(Proposition~\ref{l2.102}), but also show some remarkable
differences (Proposition~\ref{p2.5} and Remark~\ref{different}). This
result is surprising and should be taken as a warning that these pairs
of functors are quite different.

For a right/left exact endofunctor $F$ on $\Oo_0$ we denote by
$\cL F$/$\cR F$ its derived functor with $i$-th (co)homology
$\cL_i F$/$\cR^i F$.

\subsection{An action of the singular braid monoid}\label{s1.2}

We would like to start the description of our results with the very
surprising fact that some of these functors give rise to a categorification
of the singular braid monoid. In reality, this is an application of the technique,
developed during the proofs of the other statements. Baez (\cite{Baez}) and Birman
(\cite{Birman}) defined the so-called singular braid monoid with generators
\begin{displaymath}
\xymatrix@!=0.1pc{
&\text{\tiny 1}\ar@{-}[dd]& & \text{\tiny i-1}\ar@{-}[dd] & \text{\tiny i}\ar@{-}[ddrr]&&
\text{\tiny i+1}\ar@{-}'[dl][ddll] &  \text{\tiny i+2}\ar@{-}[dd] & &
\text{\tiny n}\ar@{-}[dd]\\
\sigma_{\text{\tiny i}}=&&\dots&&&&&&\dots&\\
&\text{\tiny 1}&&\text{\tiny i-1}&\text{\tiny i}&&\text{\tiny i+1}
&\text{\tiny i+2}&&\text{\tiny n}
}
\end{displaymath}
\begin{displaymath}
\xymatrix@!=0.1pc{
&\text{\tiny 1}\ar@{-}[dd]& & \text{\tiny i-1}\ar@{-}[dd] & \text{\tiny i}\ar@{-}'[dr][ddrr]&&
\text{\tiny i+1}\ar@{-}[ddll] &  \text{\tiny i+2}\ar@{-}[dd] & &
\text{\tiny n}\ar@{-}[dd]\\
\sigma_{\text{\tiny i}}^{-1}=&&\dots&&&&&&\dots&\\
&\text{\tiny 1}&&\text{\tiny i-1}&\text{\tiny i}&&\text{\tiny i+1}
&\text{\tiny i+2}&&\text{\tiny n}
}
\end{displaymath}
\begin{displaymath}
\xymatrix@!=0.1pc{
&\text{\tiny 1}\ar@{-}[dd]& & \text{\tiny i-1}\ar@{-}[dd] & \text{\tiny i}\ar@{-}[ddrr]&&
\text{\tiny i+1}\ar@{-}[ddll] &  \text{\tiny i+2}\ar@{-}[dd] & &
\text{\tiny n}\ar@{-}[dd]\\
\tau_{\text{\tiny i}}=&&\dots&&&\bullet&&&\dots&\\
&\text{\tiny 1}&&\text{\tiny i-1}&\text{\tiny i}&&\text{\tiny i+1}
&\text{\tiny i+2}&&\text{\tiny n}
}
\end{displaymath}
(for the definition see Section~\ref{BraidMonoid}) which has connections to
Vassiliev-invariants (see e.g. \cite{Birman}, \cite{GP}, \cite{Ve})
and for which the word problem is solved (e.g. \cite{Corran}, \cite{Orekov}).
The following result
suggests a strong link between our functors and invariants of knots, in the
spirit of \cite{BFK}, i.e. we expect that the following result is the first
step in defining a functorial version of Vassiliev invariants.

\begin{theorem}
\label{singularbraid}
Let $\mg=\SL_n$ and let $s_i$ ($1\leq i\leq n-1$) be the set of simple
reflections in the usual ordering. We consider the corresponding functors
$\CC_{s_i}$, $\KK_{s_i}$ and $\theta_{s_i}$. The set of functors
$\{\cL\CC_{s_i}, \cR\KK_{s_i}, \theta_{s_i}\}$ defines a (weak) action of
the singular braid monoid, with generators $\{\sigma_i, \sigma_i^{-1},
\tau_i\}$, on the bounded derived category of $\cO_0$.
\end{theorem}

\subsection{The interplay between the functors}\label{s1.3}

\begin{theorem}\label{tmain2}
There are the following isomorphisms of functors:
\begin{enumerate}
\item\label{tmain2.105}  $\mathcal{R}^1\KK\cong \wZZ$.
\item\label{tmain2.1} $\mathcal{R}^1\GG\cong \ZZ$, in particular
$\mathcal{R}^1\GG\cong \ID$ on $\Oo_0^s$.
\item\label{tmain2.2} $\mathcal{L}_1\ZZ\cong \FF$,
in particular $\FF\cong\FF'$.
\item\label{tmain2.3}
$\mathcal{R}^i\GG^2\cong\begin{cases}
\ZZ\GG&\text{if $i=1$,}\\
\ZZ&\text{if $i=2$,}\\
0&\text{if $i>2$}.
\end{cases}
\quad\quad\text{and}\quad\quad
\mathcal{R}^i\KK^2\cong\begin{cases}
\wZZ\KK&\text{if $i=1$,}\\
\wZZ&\text{if $i=2$,}\\
0&\text{if $i>2$}.
\end{cases}
$
\end{enumerate}
Dual statements hold for $\ZZ'$, $\TT$, $\hat{\ZZ}'$, and $\CC$.
\end{theorem}

\begin{remark}
\label{derived1}
{\rm
$\mathcal{R}^i\GG\cong 0$ for $i>1$ by \cite[Theorem~2.2 and Theorem~4.1]{AS};
$\mathcal{L}^2\ZZ\cong \ZZ'$ and $\mathcal{L}^i\ZZ\cong 0$ if $i>2$
follows from \cite[Corollary~5.2]{EW}, and $\mathcal{R}^i\KK\cong 0$ for $i>1$
follows from  \cite[Lemma~5.2 and Proposition~5.3]{MS}.
\erem}
\end{remark}

\begin{remark}
\label{different}
{\rm
The derived functor $\mathcal{L}\wZZ$ has a more complicated structure
than $\mathcal{L}\ZZ$. This is already evident for the Lie algebra
$\mathfrak{sl}_3$. In fact, by a direct calculation one can show that
in this case $\mathcal{L}_{6}\wZZ\neq 0$. It follows that, in general,
there is no involutive exact equivalence $F$ on $\Oo_0$ sending $L(x)$ to
$L(x^{-1})$. The same statement can also be obtained using the
following general argument:

Let $A$ be a finite-dimensional  associative algebra and $\Lambda$ be
an indexing set of the isoclasses $S(\lambda)$, $\lambda\in\Lambda$
of simple $A$-modules. Assume that $F$ is an
exact equivalence on $A\mathrm{-mod}$ such that $F(S(\lambda))\cong
S(\sigma(\lambda))$ for some permutation $\sigma$ on $\Lambda$. For
$J\subset \Lambda$ let $Z_J$ denote the functor given by taking the maximal
quotient containing only simple subquotients indexed by $J$. Then
it is easy to see that the functors $F^{-1}Z_{\sigma(J)}F$ and
$Z_J$ are isomorphic.

Let $\mathfrak{g}=\mathfrak{sl}_3$ and $s$, $t$ be the two simple reflections. Let $J=\{e,t,ts\}$, $\hat{J}=\{e,t,st\}$ and
$J'=\{e,s,ts\}$. Then
$J\cong\hat{J}$ via $w\mapsto w^{-1}$ and $J\cong J'$ via $ww_0\mapsto
w^{-1}w_0$. By definition we have $\ZZ= Z_J$,
$\hat{\ZZ}=Z_{\hat{J}}$, and $\hat{\ZZ}_t=Z_{J'}$. It is easy to check that
$\ZZ P(t)$ has length $4$, but both,  $\hat{\ZZ} P(t^{-1})$ and $\hat{\ZZ}_t
P(s)=\hat{\ZZ}_t P((st)^{-1}w_0)$, have length $3$. In particular, there is
neither an involutive exact equivalence sending $L(x)$ to $L(x^{-1})$, nor an
involutive exact equivalence sending $L(xw_0)$ to $L(x^{-1}w_0)$. From the
point of view of Kazhdan-Lusztig combinatorics this could not have been expected. In particular, it shows that the corresponding statements in the
literature (e.g. \cite[Existence of $\varepsilon$ in Section 4.3]{JCompl}
and the existence of $T$ claimed in \cite[Lemma 6]{SoInv}) are not correct.
In fact the construction proposed in the proof of \cite[Lemma 6]{SoInv}
gives the identity functor.
\erem}
\end{remark}

We describe the monoids generated by $\{\GG,\TT\}$ and $\{\CC,\KK\}$
respectively. Recall (see e.g. \cite[Chapter~II]{La}), that for a monoid,
$S$, and for $x,y\in S$, we have $x{\bf R}y$ or $x{\bf L}y$ if and only if
$xS=yS$ and $Sx=Sy$ respectively.

\begin{theorem}\label{tmain3}
The functors $\TT$ and $\GG$ satisfy the relations
\begin{gather*}
\TT\GG\TT\cong\TT,\; \GG\TT\GG\cong\GG,\;
\TT^3\cong\TT^2,\;\GG^3\cong\GG^2,\\
\TT^2\GG\cong\TT^2,\; \GG^2\TT\cong\GG^2,\;
\TT\GG^2\cong \GG\TT^2,
\end{gather*}
and their isoclasses generate the monoid
$\mathcal{S}=\{\ID,\TT,\GG,\TT\GG,\GG\TT,\TT^2,\GG^2,\TT\GG^2\}$ of
(isoclasses of) functors.
The columns and rows of the following egg-box diagrams
represent respectively  Green's relations ${\bf R}$ and ${\bf L}$,
on $\mathcal{S}$:
\begin{displaymath}
\begin{array}{|c|}
\hline
\ID\\
\hline
\end{array}\quad\quad
\begin{array}{|c|c|}
\hline
\GG&\TT\GG\\
\hline
\GG\TT&\TT\\
\hline
\end{array}\quad\quad
\begin{array}{|c|c|c|}
\hline
\GG^2&\TT^2&\GG\TT^2\\
\hline
\end{array}
\end{displaymath}
\end{theorem}

\begin{theorem}\label{tmain3new}
The functors $\CC$ and $\KK$ satisfy the relations
\begin{gather*}
\CC\KK\CC\cong\KK,\;\KK\CC\KK\cong\KK,\;
\CC^3\KK\cong\CC^2,\;
\KK^3\CC\cong\KK^2,\;\\
\CC^2\KK^2\CC\cong\CC^2\KK,\;\KK^2\CC^2\KK\cong\KK^2\CC,\;
\CC\KK^2\CC^2\cong\KK\CC^2,\; \KK\CC^2\KK^2\cong\CC\KK^2.
\end{gather*}
Assume that $s$ does not correspond to an $\mathfrak{sl}_2$-direct
summand of $\mathfrak{g}$. Then the
isoclasses of the functors $\CC$ and $\KK$ generate the
(infinite) monoid
\begin{equation*}
\mathcal{\hat{S}}=\{\ID, \KK\CC^2\KK\cong\CC\KK^2\CC, \KK^i, \CC^i,
\KK\CC^i,\CC\KK^i, \KK^2\CC^i,\CC^2\KK^i:i>0\}.
\end{equation*}
The columns and rows of the following egg-box diagrams
represent respectively  Green's relations ${\bf R}$ and ${\bf L}$,
on $\mathcal{\hat{S}}$:
\begin{displaymath}
\begin{array}{|c|}
\hline
\ID\\
\hline
\end{array}\quad\quad
\begin{array}{|c|c|}
\hline
\KK&\CC\KK\\
\hline
\KK\CC&\CC\\
\hline
\end{array}\quad\quad
\begin{array}{|c|c|cc|}
\hline
\CC^i,i>1, & \KK^i,i>1, & \CC\KK^i, & \KK\CC^i, \\
\CC^2\KK^i,i>0 & \KK^2\CC^i,i>0& i>1,&
\KK\CC^2\KK\\
\hline
\end{array}
\end{displaymath}
The only idempotents in $\mathcal{\hat{S}}$ are
$\ID$, $\KK\CC$, $\CC\KK$, $\CC^2\KK^2$, $\KK^2\CC^2$, $\KK\CC^2\KK$.
\end{theorem}

\begin{remark}\label{serre}
{\rm
If $\mathfrak{g}=\mathfrak{sl}_2$, then one can show that
$\cL\CC^2$ is the Serre functor, which implies  $\CC^2\cong \CC^4$,
see \cite{KMS}.
}
\end{remark}

Before describing morphism spaces between such functors, we
want to give an impression of their rather complex interplay.
We need some preparations to formulate the corresponding
Theorem~\ref{tmain1}, in which we show relations between functors
from $\mathcal{S}$.

According to  \cite[Remark~5.7]{AS}, $\TT$ is left adjoint to
$\GG$ and $\rt$ is up to a scalar the composition of $\TT(\rg)$
with the adjunction morphism $\TT\GG{\longrightarrow}\ID$. We fix
$\raa'\in \Hom(\TT\GG,\ID)$ such that $\rt=\raa'\circ \TT(\rg)$ and
set $\raa=\rd(\raa')_{\rd}$ (the existence of $\raa'$ also
follows from the independent result $\Hom(\TT\GG,\ID)\cong \Coinv$ of
Theorem~\ref{tmain4} which ensures that up to a scalar there is only
one natural transformation ``of degree zero''). Let $\rz:\ID\surj\ZZ$,
and $\rp:\GG\surj\FF$ be the natural projections,
$\ri=\rd(\rp)_{\rd}$,
$\rk'=(\TT^2(\rg))^{-1}\circ\ri_{\TT\GG}$, and
$\rk=\rd(\rk')_{\rd}$. We will see later that all these maps are
well-defined.

\begin{theorem}\label{tmain1}
Figure~\ref{P2} presents a diagram of endofunctors on $\Oo_0$,
where all the sequences labeled by numbers are exact.
\begin{figure}[htbp]
  \centering $\hspace{-1cm}
\xymatrix@!=0.2pc{
 &&\FF'\TT\ar@{^{(}->}[lldd]_{\ZZ'(\ri_{\TT})}|{{}\text{\bf 5}}
\ar@{^{(}->}[rrdd]^{\ri_{\TT}}|{{}\text{\bf 7}}
\ar@{^{(}->}[dddd]_<<<<<<<<<<<<<<<<<<<<<<<<<<<{\FF'\TT(\rg)}|<<<<<<{\text{\bf 6}}
&& &&{\ID}
\ar[ddrr]^<<<<{\GG(\rg)\circ\rg}
&& && \FF\GG
\ar@/_7pc/@{^{(}.>}[dddddddddlllllll]_>>>{{\ri_{\GG}}}|>>>>>>>>>>>>>>>{{}\text{\bf 10}}&&\\
&&&&&&&&&&&&\\
\ZZ'\TT^2 \ar@{->>}[dddd]_>>>>>>>>>{\ZZ'\TT(\rt)}|{{}\text{\bf 5}}
\ar@{^{(}->}[rrrr]^>>>>>>>>>>{\rz'_{\TT^2}}|<<<<<<<<{{}\text{\bf 1}}
&& && \TT^2
\ar[rrrr]^{\GG(\rg)\circ\rg\circ\rt\circ\TT(\rt)}|{{}\text{\bf 1}}
\ar[rruu]^>>>>{\rt\circ\TT(\rt)}\ar@{->>}[dddd]_>>>>{\rt_{\TT}}|<<<<<<{{}\text{\bf 7}}
&& && \GG^2
\ar@{->>}[rrrr]^>>>>>>>>>>{\rz_{\GG^2}}|<<<<<<<{{}\text{\bf 1}}
\ar@{->>}[rruu]^{\rp_{\GG}}|{{}\text{\bf 7$'$}}\ar@{->>}[rrdd]^{\rk}
&& && \ZZ\GG^2\ar@{->>}[lluu]_{\ZZ(\rp_{\GG})}|{{}\text{\bf 5$'$}} \\
&&&&&&&&&&&&\\
 && \FF'\TT\GG \ar@{^{(}->}[rruu]^{\rk'}
\ar@{->>}[rrr]^>>>>>>>>>>>>>>>>>{\FF'(\raa')}|{{}\text{\bf 6}}
&& & \FF'\ar@{^{(}->}[ddl]^{\ri}|<<<<<<{{}\text{\bf 8 3$'$}}
\ar@{.>}[rr]_{\mathrm\alpha}^{\sim}
 && \FF
\ar@{^{(}->}[rrr]^<<<<<<<<<<<<<<<<<{\FF(\raa)}|{{}\text{\bf 6$'$}}& &&
\FF\GG\TT
\ar@{->>}[uuuu]_>>>>>>>>>>>>>>>>>>>>>>>>>>>{\FF\GG(\rt)}|>>>>>>{{}\text{\bf 6$'$}}   && \\
&&&&&&&&&&&&\\
\ZZ'\TT\ar@{^{(}->}[rrrr]^<<<<<<<<<{\ZZ'_{\TT}}|{{}{\text{{\bf 2 9$'$}}}}
\ar[dddd]_{\ZZ'(\rt)}
 && && \TT \ar[rrrr]^<<<<<<<<<<<{\rg\circ\rt}|{{}\text{\bf 2}}
\ar[ddddlll]_{\rt}|<<<<<<{{}\text{\bf 3$'$}}
\ar@{->>}[lddd]^{\TT(\rg)}|<<<<<<<<{{}\text{\bf 8}}
\ar@{->>}[ddr]^{\TT(\rg_{\GG}\circ\rg)}|{{}\text{\bf 9$'$}}
&& && \GG
\ar@{->>}[rrrr]^<<<<<<<<<<<<<<<{\rz_{\GG}}|{{}{\text{\bf 2 9}}}
\ar@{_{(}->}[uuuu]^>>>>{\rg_{\GG}}|<<<<<<{{}\text{\bf
 7$'$}}\ar@{->>}[luu]^{\rp}|>>>>>>{{}{\text{\bf 3 8$'$}}}
&& && \ZZ\GG
\ar@{^{(}->}[uuuu]^<<<<<<<<<{\ZZ\GG(\rg)}|{{}\text{\bf 5$'$}}\\
&&&&&&&&&&&&\\
 &&  & & &
\TT\GG^2\ar[rr]_{\rh}^{\sim}
& &
\GG\TT^2\ar@{^{(}->}[ruu]^{\GG(\rt\circ\rt_{\TT})}|{{}\text{\bf 9}}
\ar@{^{(}->}[rrd]_{\GG(\rt_{\TT})}|<<<<<<{{}\text{\bf 10$'$}}
& &
& && \\
&&&
\TT\GG\ar@{_{(}->}[lld]^{\raa'}|<<<<<<{{}\text{\bf 4$'$}}
\ar@{->>}[rru]_{\TT(\rg_{\GG})}|{{}\text{\bf 10}}
&&&&&&
\GG\TT\ar@{_{(}->}[luuu]^{\GG(\rt)}|<<<<<<<<{{}\text{\bf 8$'$}}
\ar@/_7pc/@{.>>}[llllllluuuuuuuuu]_>>>{{\rp_{\TT}}}|>>>>>>>>>>>>>>>>>>>>>>>>>>>>>>>>>>>>>>>>>>>>>>>>>>>>>>>{{}\text{\bf 10$'$}}
&&&\\
\ZZ'\ar@{^{(}->}[r]^{\rz'}_{{}{\text{\bf 3 4 }}}
%\ar@{.>}@<-0.1em>
&
\ID \ar@{=}[rrrrrrrrrr]
&&& && && &&&
\ID
\ar@{->>}[r]^{\rz}_{{}{\text{\bf 3$'$ 4$'$}}}
\ar@{->>}[llu]^{\raa}|>>>>>>{{}{\text{\bf 4}}}
\ar[uuuulll]_{\rg}|>>>>>{{}{\text{\bf 3}}}
& \ZZ \ar[uuuu]^{\ZZ(\rg)}\\
}
$
 \caption{Commutative diagram involving $\TT$ and $\GG$}
  \label{P2}
\end{figure}
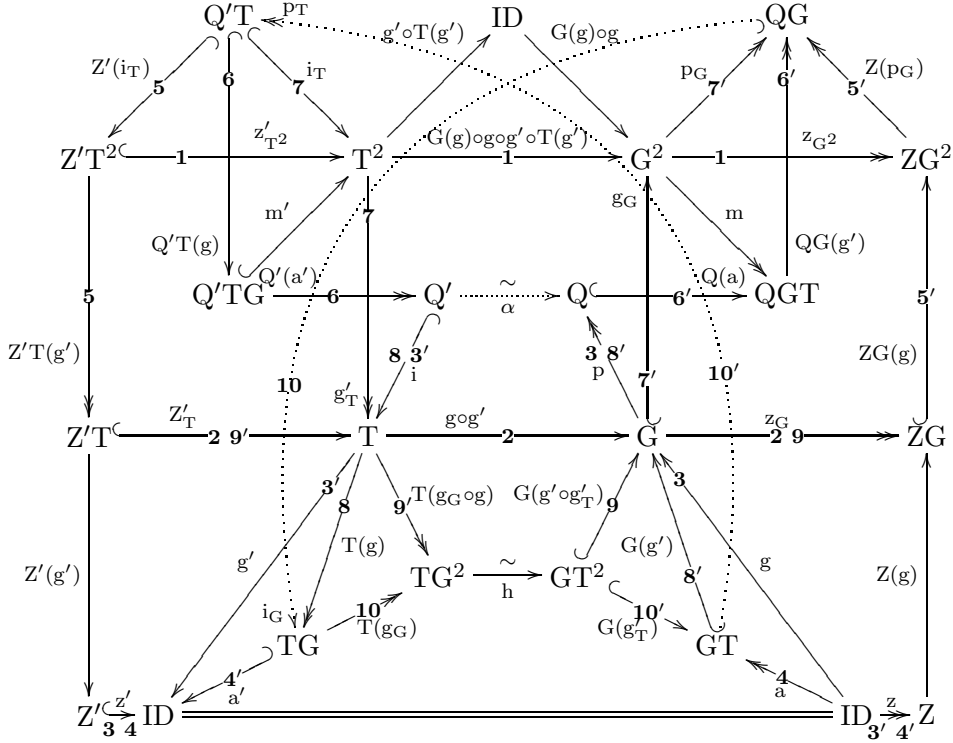
Moreover, one can choose isomorphisms ${\mathrm \alpha}$ and $\rh$ as
indicated such that all configurations containing only solid arrows commute.
\end{theorem}

\subsection{Natural transformations}\label{s1.101}

We prove the following result on natural transformations between
arbitrary compositions of $G$ and $T$:

\begin{theorem}\label{tmain4}
\begin{enumerate}
\item For $\XX\in \mathcal{S}$ there
is a ring isomorphism $\End(X)\cong \Coinv$.
\item For $\XX$, $\YY\in\mathcal{S}$ we have $\Hom(\XX,\YY)\not=0$ and
this space is given by the $\XX$-row and $\YY$-column entry in
the following table:
\begin{displaymath}
\begin{array}{|c||c|c|c|c|c|c|c|c|}
\hline
\XX\setminus\YY & \ID & \GG & \TT & \GG\TT & \TT\GG
& \GG^2 & \TT^2 & \GG\TT^2 \\
\hline\hline
\ID & \Coinv & \Coinv & 1 & \Coinv & 2 & \Coinv & 3 & 4 \\
\hline
\GG & 1 & \Coinv & 5 & 4 & 1 & \Coinv & 6 & 4 \\
\hline
\TT & \Coinv & \Coinv & \Coinv & \Coinv & \Coinv & \Coinv & 4 &
\Coinv\\
\hline
\GG\TT & 2 & \Coinv & 1 & \Coinv & 7 & \Coinv & 8 & 4 \\
\hline
\TT\GG & \Coinv & \Coinv & 4 & \Coinv & \Coinv & \Coinv & 4 &
\Coinv \\
\hline
\GG^2 & 3 & 4 & 6 & 4 & 8 & \Coinv & 9 & 4 \\
\hline
\TT^2 & \Coinv & \Coinv & \Coinv & \Coinv & \Coinv & \Coinv &
\Coinv & \Coinv \\
\hline
\GG\TT^2 & 4 & \Coinv & 4 & \Coinv & 4 & \Coinv & 4 & \Coinv \\
\hline
\end{array}.
\end{displaymath}
The spaces described by the same number are isomorphic and we have the
following inclusions:
\begin{eqnarray*}
\xymatrix{
{\bf A}:\quad 7\ar@{^{(}->}[r] & 2\ar@{^{(}->}[r] & 4\ar@{^{(}->}[r] &
\Coinv
&
{\bf B}:\quad
8\ar@{^{(}->}[r] & 3\ar@{^{(}->}[r]\ar@{^{(}->}[rd] & 6 \\
& 1\ar@{^{(}->}[ru] &  &
 & & & 9
\\
{\bf C}:\quad
\END(\ID_{\cO_0^s})\ar@{^{(}->}[r]& 5}
\end{eqnarray*}
\item There is an isomorphism of rings $\End(\ZZ)\cong \END(\ID_{\cO_0^s})$.
\end{enumerate}
\end{theorem}

We describe the endomorphism spaces of the elements from
$\mathcal{\hat{S}}$ and natural transformations between the
idempotents in the following theorem:

\begin{theorem}\label{tmain4.101}
\begin{enumerate}
\item For $\XX\in\mathcal{\hat{S}}$ there
is a ring isomorphism $\End(X)\cong \Coinv$.
\item For idempotents $\XX$, $\YY\in\mathcal{\hat{S}}$
the space $\Hom(\XX,\YY)$ is given by the $\XX$-row and
$\YY$-column entry in the following table:
\begin{displaymath}
\begin{array}{|c||c|c|c|c|c|c|}
\hline
\XX\setminus\YY & \ID & \CC\KK & \KK\CC & \CC^2\KK^2 & \KK^2\CC^2
& \KK\CC^2\KK \\
\hline\hline
\ID & \Coinv & 1 & \Coinv & 2 & \Coinv & 3 \\
\hline
\CC\KK & \Coinv & \Coinv & \Coinv & 4 & \Coinv & \Coinv \\
\hline
\KK\CC & 1 & 5 & \Coinv & 2 & \Coinv & 3 \\
\hline
\CC^2\KK^2 & \Coinv & \Coinv & \Coinv & \Coinv & \Coinv & \Coinv \\
\hline
\KK^2\CC^2 & 2 & 2 & 4 & 6 & \Coinv & 4 \\
\hline
\KK\CC^2\KK & 3 & 3 & \Coinv & 4 & \Coinv & \Coinv \\
\hline
\end{array}.
\end{displaymath}
The spaces described by the same number are isomorphic and we have the
following inclusions:
\begin{displaymath}
5\hookrightarrow 1\hookrightarrow 3\hookrightarrow \Coinv, \quad\quad
4\hookrightarrow \Coinv.
\end{displaymath}
\end{enumerate}
\end{theorem}

\begin{remark}{\rm
The coinvariant algebra has a natural $\mZ$-grading given by putting
$\mh$ in degree one. Using the graded versions of $\CC$ and $\KK$
from \cite[7.1]{MS} (and a similar construction for $\GG$ and $\TT$) we get
isomorphisms of graded vector spaces as listed in the theorem.
\erem}
\end{remark}

\subsection{Tilting modules}\label{s1.4}

Let $\mathcal{P}=\oplus_{x\in W} P(x)$ be a minimal projective
generator of $\cO_0$ and set $\mathcal{I}=\rd\mathcal{P}$. For
$M\in\cO_0$ the category $\mathrm{Add}(M)$ is defined as the full
subcategory of $\cO_0$, which consists of all direct summands of
all finite direct sums of copies of $M$.  Recall (see \cite[Introduction]{Wa}) that
$M\in\cO_0$ is called a {\em generalized tilting module} if
$\Ext^{>0}_{\cO_0}(M,M)=0$ and if $\mathcal{P}$ has a finite
$\mathrm{Add}(M)$-coresolution, i.e. there  exists an exact
sequence $0\to \mathcal{P}\to M_0\to\dots\to M_k\to 0$ of finite
length $k$ with $M_i\in\mathrm{Add}(M)$ for $1\leq i\leq k$. If,
additionally, the projective dimension of $M$ is one then $M$ is
called a {\em classical tilting module}, see \cite[page 399]{HR}. Dual notions
define generalized and classical cotilting modules. For a fixed
reduced expression $w=s_1\cdots s_k\in W$ we set $\TT_w=\TT_{s_1}
\cdots \TT_{s_k}$ and $\GG_w=\GG_{s_1}\cdots \GG_{s_k}$. The resulting
functors are (up to isomorphism) independent of the chosen reduced
expression (see \cite[2.9]{JCompl}, \cite[Section~6]{KM}). The following result describes
a lattice of (generalized) tilting and cotilting modules in $\Oo_0$
constructed using twisting and completion functors.

\begin{theorem}\label{tmain5}
Let $w\in W$.
\begin{enumerate}
\item Each of the modules $\mathcal{P}^{w}=\TT_w \mathcal{P}$ and
$\mathcal{I}^w=\GG_w\mathcal{I}$ is both, a generalized tilting module
and a generalized cotilting module.
\item We have the following  equalities for
projective and injective dimensions:
$\mathrm{projdim}
(\mathcal{P}^w)=\mathrm{injdim}(\mathcal{I}^w)=l(w)$ and
$\mathrm{injdim}(\mathcal{P}^w)=\mathrm{projdim}(\mathcal{I}^w)=
2l(w_0)-l(w)$.
In particular, if $s$ is a simple reflection then $\mathcal{P}^s$
($\mathcal{I}^s$ resp.) is a classical (co)tilting module.
\item $\TT_{w}\mathcal{P}^{w_0}\cong\mathcal{I}^{ww_0}$ and
$\GG_{w}\mathcal{I}^{w_0}\cong\mathcal{P}^{ww_0}$. In particular,
$\mathcal{P}^{w_0}\cong \mathcal{I}^{w_0}\cong \mathcal{T}$ is
the characteristic (co)tilting module in $\Oo_0$.
\end{enumerate}
\end{theorem}

\begin{remark}{\rm Let $x\in W$ be fixed.
The module $\TT_{x}\TT_{w_0}\mathcal{P}\cong \TT_{x}\mathcal{P}^{w_0}
\cong \TT_{x}\mathcal{T}$ is the direct sum of all
{\em $x$-twisted tilting modules}
as defined in \cite[Section~5]{St4} and characterized by certain
vanishing conditions with respect to twisted Verma modules. If $x=e$
we get the sum of all (usual) tilting modules. The twisting functors
define naturally maps as follows:
\begin{multline*}
\{\text{indec. projectives}\}
\overset{\TT_x}{\longrightarrow}
\{x\text{-twisted indec. projectives}\}
\overset{\TT_{w_0x^{-1}}}{\longrightarrow}\\
\overset{\TT_{w_0x^{-1}}}{\longrightarrow}
\{\text{($e$-twisted) tiltings}\}
\overset{\TT_x}{\longrightarrow}
\{x\text{-twisted tiltings}\} =\\
=\{xw_0\text{-completed indec. injectives}\}
\overset{T_{w_0x^{-1}}}{\longrightarrow}
\{\text{indec. injectives}\}.
\end{multline*}
The maps are all bijections, their inverses induced by the
corresponding completion functors.
\erem}
\end{remark}

For a reduced expression  $w=s_ks_{k-1}\cdots s_1\in W$ we
set $\CC_w=\CC_{s_1}\cdots \CC_{s_k}$ and
$\KK_w=\KK_{s_1}\cdots\KK_{s_k}$. Up to isomorphism, the functors do
not depend on the chosen reduced expression, see \cite[Lemma~5.10]{MS}. We will
prove the following analog of the previous theorem:

\begin{theorem}\label{tmain6}
Let $w\in W$.
\begin{enumerate}
\item\label{tmain6.1} Each of the modules
${}_{}^{w}\mathcal{P}=\CC_w \mathcal{P}$ and
${}_{}^w\mathcal{I}=\KK_w\mathcal{I}$ is both, a generalized tilting
module and a generalized cotilting module.
\item\label{tmain6.2} We have the following
equalities for projective and injective dimensions:
$\mathrm{projdim}({}_{}^w\mathcal{P})=
\mathrm{injdim}({}_{}^w\mathcal{I})=l(w)$ and
$\mathrm{injdim}({}_{}^{w}\mathcal{P})=
\mathrm{projdim}({}_{}^{w}\mathcal{I})= 2l(w_0)-l(w)$.
In particular, ${}_{}^s\mathcal{P}$ (and ${}_{}^s\mathcal{I}$ resp.)
is a classical (co-)tilting module for any simple reflection
$s\in W$.
\item\label{tmain6.3}
$\CC_{w}({}_{}^{w_0}\mathcal{P})\cong{}_{}^{w^{-1}w_0}\mathcal{I}$
and
$\KK_{w}({}_{}^{w_0}\mathcal{I})\cong{}_{}^{w^{-1}w_0}\mathcal{P}$.
In particular, ${}_{}^{w_0}\mathcal{P}\cong {}_{}^{w_0}\mathcal{I}\cong
\mathcal{T}$ is the characteristic (co)tilting module in $\Oo_0$.
\end{enumerate}
\end{theorem}

\begin{question}{\rm
According to \cite[Theorem~5.4]{AR} every generalized tilting module $T$ for an
associative algebra $A$ corresponds to a resolving and contravariantly
finite subcategory in $A\mathrm{-mod}$ consisting of all $A$-modules
admitting a finite coresolution by $\mathrm{Add}(T)$. What are the
subcategories of $\Oo_0$, which correspond to the various generalized
tilting objects from above?
}
\end{question}

\section{Preliminary properties of our functors}\label{s2}

In this section we collect some fundamental statements concerning
natural transformations between our functors. As a corollary we get
a short argument for the existence of an isomorphism $T\cong\GG'$
(which was originally proved in \cite[Theorem~4]{KM}).

By \cite[Endomorphismensatz]{S} we have $\END_\mg(P(w_0))\cong \Coinv$, and thus we can
define the functor $\mathbb{V}:\Oo_0\to \Coinv\mathrm{-mod}$,
$M\mapsto \Hom_{\g}(P(w_0),M)$. Let $\tilde\GG$ denote the right-adjoint of
$\TT$, which exists by \cite[Section~4]{AS}.

\begin{lemma}\label{l2.1}
$\mathbb{V}\tilde\GG\cong \mathbb{V}$ and $\tilde\GG\cong\ID$
when restricted to projectives.
\end{lemma}

\begin{proof}
Note that $\TT P(w_0)\cong P(w_0)$ and $\END_\mg(P(w_0))$ is given by
the action of the center $\cZ$ of  the universal enveloping algebra of
$\mg$ (\cite[Endomorphismensatz]{S}). On the other hand, the action of $\cZ$ commutes
naturally with $\TT$ by definition. This allows us to fix a natural
isomorphism $\TT\cong\ID$ on $\mathrm{Add}(P(w_0))$.
This ensures that (for any $M\in\cO_0$) the following
isomorphisms are even morphisms of $\Coinv$-modules:
\begin{eqnarray*}
\mathbb{V} M&=&
\Hom_{\g}(P(w_0),M)\cong
\Hom_{\g}(\TT P(w_0),M)\cong
\Hom_{\g}(P(w_0),\tilde\GG M)\\
&=&\mathbb{V}\tilde\GG M.
\end{eqnarray*}
All the isomorphisms are natural and the first statement follows.
Let $\tilde\mV$ denote the right-adjoint of $\mV$. By
\cite[Proposition 6]{S}  we have $\tilde{\mathbb{V}}\mathbb{V}\cong\ID$
on projectives and therefore also $\tilde\GG\cong
\tilde{\mathbb{V}}\mathbb{V}\tilde\GG\cong
\tilde{\mathbb{V}}\mathbb{V}\cong\ID$, since $\tilde\GG$ preserves the
category of projectives.
\end{proof}

We fix an isomorphism of functors $\varphi:\ID\cong \tilde\GG$ defined
on the category of projectives. For $M\in\Oo_0$ we choose a projective
presentation
\begin{displaymath}
P_1\overset{\gamma'}{\longrightarrow}
P_0\overset{\gamma}{\tto} M.
\end{displaymath}
Then the left square of the following
diagram commutes and induces the map $\varphi_M$ as indicated:
\begin{displaymath}
\xymatrix{
\tilde\GG P_1 \ar[rr]^{\tilde\GG \gamma'}&&
\tilde\GG P_0\ar[rr]^{\tilde\GG \gamma}&& \tilde\GG M\\
P_1\ar[rr]^{\gamma'}\ar[u]^{\varphi_{P_1}}&&
P_0\ar@{->>}[rr]^{\gamma}\ar[u]^{\varphi_{P_0}}&&
M\ar@{-->}[u]^{\varphi_M}
}.
\end{displaymath}

\begin{lemma}\label{l2.3}
The maps $\phi_M$, $M\in \Oo_0$, define a natural transformation
from $\ID$ to $\tilde\GG$.
\end{lemma}

\begin{proof}
First we have to check that $\phi_M$ is independent of the chosen
presentation. Let $Q_1\overset{\beta'}{\longrightarrow}
Q_0\overset{\beta}{\tto} M$ be another projective presentation of
$M$. Consider the commutative diagram:
\begin{displaymath}
\xymatrix{
\tilde\GG P_1 \ar[rr]^{\tilde\GG \gamma'}&&
\tilde\GG P_0\ar[rr]^{\tilde\GG \gamma}&& \tilde\GG M\\
P_1\ar[rr]^{\gamma'}\ar[u]^{\varphi_{P_1}}&&
P_0\ar@{->>}[rr]^{\gamma}\ar[u]^{\varphi_{P_0}}&&
M\ar@{..>}[u]^{h}\ar@{=}[d] \\
Q_1\ar[rr]^{\beta'}\ar[d]^{\varphi_{Q_1}}\ar@{..>}[u]^{\xi'}&&
Q_0\ar@{->>}[rr]^{\beta}\ar[d]^{\varphi_{Q_0}}\ar@{..>}[u]^{\xi}&&
M\ar@{..>}[d]^{h'} \\
\tilde\GG Q_1 \ar[rr]^{\tilde\GG \beta'}&&
\tilde\GG Q_0\ar[rr]^{\tilde\GG \beta}&& \tilde\GG M\\
},
\end{displaymath}
where the projectivity of $Q_1$ and $Q_0$ is used to get $\xi'$ and
$\xi$ such that the diagram is commutative. Since $\xi$ is a map
between projectives, we obtain $\tilde\GG\xi\circ\varphi_{Q_0}=
\varphi_{P_0}\circ \xi$. Hence
\begin{displaymath}
h'\circ\beta=\tilde\GG\beta\circ \phi_{Q_0}=\tilde\GG\gamma\circ
\tilde\GG\xi\circ\phi_{Q_0}=\tilde\GG\gamma\circ
\phi_{P_0}\circ\xi=h\circ \gamma\circ\xi=h\circ  \beta,
\end{displaymath}
by the commutativity of the diagram. Since $\beta$ is surjective, we
obtain $h=h'$. Hence, $\phi_M$ is well-defined. The naturality
follows by standard arguments.
\end{proof}

\begin{proposition}\label{p2.4}
$\GG$ is right adjoint to $\TT$. In particular, there exists a natural
transformation $\TT\rightarrow\ID$ non-vanishing on Verma modules.
\end{proposition}

\begin{proof}
Lemma~\ref{l2.3} implies the existence of a non-trivial natural
transformation $\TT\to\ID$ as assumed in
\cite[Proposition~5.4]{AS}. The statement now follows from
\cite[Proposition~5.4]{AS} and \cite[Lemma~1]{KM}.
\end{proof}

\begin{proposition}\label{l2.102}
\begin{enumerate}[(1)]
\item \label{1} $(\TT,\GG)$ is an adjoint pair of functors. The
adjunction morphism $\ADJ_{\TT}:\TT\GG\to \ID$ is injective with
cokernel $\ZZ$, and the adjunction morphism $\ADJ^{\TT}:\ID\to
\GG\TT$ is surjective with kernel $\ZZ'$.
\item \label{2} $(\CC,\KK)$ is an adjoint pair of functors. The
adjunction morphism $\ADJ_{\CC}:\CC\KK\to \ID$ is injective with
cokernel $\hat{\ZZ}$, and the adjunction morphism
$\ADJ^{\CC}:\ID\to \KK\CC$ is surjective with kernel $\hat{\ZZ'}$.
\item \label{3} The functors $\TT\GG$ and $\GG\TT$ preserve both
surjections and injections (but are neither left nor right exact).
\item \label{4} The functors $\CC\KK$ and $\KK\CC$ preserve both
surjections and injections (but are neither left nor right exact).
\end{enumerate}
\end{proposition}

\begin{remark}{\rm
The twisting functor $\TT$ can be described and generalized as
follows (this was also observed by W. Soergel): We consider $\cO_0$ as
the category $\MOF-A$ of finitely generated right modules over
$A=\END_\mg(\mathcal{P})$ with endofunctor $\TT$. To each simple object
$L(w)$ we have the corresponding primitive idempotent $e_w\in A$. Let
$e$ be the sum of all $e_w$ taken over all $w$ such that $T L(w)\not=0$
(see \cite[Proposition~5.1]{AS}) and define $\tilde{\TT}={}_-\otimes_A
AeA:\mathrm{mod-}A\rightarrow\mathrm{mod-}A$.
Using the  definitions and \cite[Proposition~5.3, Theorem~2.2 and
Corollary~5.2]{AS} we get $\TT(A_{A})\cong \tilde{\TT}(A_{A})$ and the
inclusion $AeA\hookrightarrow A$ induces a non-trivial element
$\phi\in\Hom(\tilde{\TT},\ID)$. Applying \cite[Lemma~1]{KM} one gets
$\tilde{\TT}\cong \TT$ as endofunctors of $\MOF-A$. This description
allows a generalization of twisting functors to a very general setting.
The definitions immediately show that the cokernel of $\phi_M$ is
always the largest quotient of $M$, such that $\HOM_A(eA,M)=0$ and
one easily derives $\tilde{\TT}^3\cong\tilde{\TT}^2$. However, if
$\tilde{\GG}$ denotes the right adjoint of $\tilde\TT$, then
the adjunction morphism $\tilde{\TT}\tilde{\GG}\rightarrow\ID$ does
not need to be injective in general.
\erem}
\end{remark}

\begin{proof}[Proof of Proposition~\ref{l2.102}]
In this proof for $M\in\Oo_0$ we denote by $[M]$ the class of $M$ in
the Grothendieck group of $\Oo_0$.

The first part is proved in \cite[Section~5]{AS}. For the
part~\eqref{3}  it is enough to show that both, $\TT\GG$ and $\GG\TT$,
preserve surjections. Assume $f\in\HOM(M,N)$ for some $M$, $N\in\cO_0$
is surjective. The adjunction morphism $\operatorname{adj}^{\TT}$ is
surjective. Then $\op{adj}_N^{\TT}\circ f=\GG\TT(f)\circ
\op{adj}_M^{\TT}$ is surjective; in particular, so is $\GG\TT(f)$.

Let $\IM$ be the image of $\GG(f)$. Then $\TT(\GG(f)):\TT\GG M\surj
\TT(\IM)$ is surjective and so is $\TT(i):\TT(\IM)\surj\TT\GG N$,
since the cokernel of $i:\IM\inj G N$ is annihilated by $\TT$. The
composition of both surjections is exactly $\TT\GG(f)$ and so we are
done: part~\eqref{3} follows.

Concerning statement~\eqref{4}, it is enough to prove the claim for
$\CC\KK$. Let us first show that $\CC\KK$ preserves inclusions. Let
$M\overset{f}{\hookrightarrow}N\overset{g}{\tto}L$ be a short exact
sequence in $\Oo_0$. Applying $\KK$  gives an
exact sequence ${\bf S}$ of the form $\KK M\hookrightarrow\KK
N\tto L'$ where $L'$ is a submodule of $\KK L$. By definition of
$\KK$, the socle of $\KK L$, and hence also of $L'$, contains only
simple modules not annihilated by $\theta_s$, hence $\cL_1\CC (L')=0$
by \cite[Section~5]{MS}. In particular, $\CC {\bf S}$ is exact, and
therefore $\CC\KK(f)$ is an inclusion.

On the other hand, applying $\KK$ to
$M\overset{f}{\hookrightarrow}N\overset{g}{\tto}L$ yields an exact
sequence ${\bf T}$ of the form
$\KK M\hookrightarrow \KK N\rightarrow \KK L\tto X$, where
$\KK X=\CC X=0$ by \cite[Proposition~5.3]{MS}. Applying
the right exact functor $\CC$ to ${\bf T}$ and using $\CC X=0$ we
obtain that $\CC\KK(g)$ is a surjection. This shows part~\eqref{4}.

By \cite[Section~5]{MS} the adjunction morphism defines an
isomorphism $\CC\KK\cong\ID$ when restricted to modules having a dual
Verma flag. Let $M\in\cO_0$ with injective cover $i:M\inj I$. Let
$\ADJ=\ADJ_{\CC}$ for the moment. Then $i\circ{\ADJ}_M={\ADJ}_{I}
\circ\CC\KK(i)$. The latter is injective, hence $\ADJ_M$ has to be
injective as well. Note that
$[\CC\KK(M)]=[\theta\KK(M)]-[\KK(M)]=
[\theta^2(M)]-[\theta(M)]-[\KK(M)]=[\theta(M)]-[\KK(M)]$
for any $M\in\cO_0$.
Hence $[M]-[\CC\KK(M)]=[\hat{\ZZ}(M)]$.
Dual statements hold for $\ADJ^{\CC}$. Part~\eqref{2} follows.
\end{proof}

The following result is surprising in comparison with
Proposition~\ref{p2.4} (note that the argument of Lemma~\ref{l2.1}
does not work if we replace $\tilde{\GG}$ by $\KK$ as $\KK$ does not
commute with the action of the center of $\Oo_0$).

\begin{proposition}\label{p2.5}
\begin{enumerate}
\item There is {\bf no} natural transformation $\rc:\CC\to\ID$
non-vanishing on Verma modules.
\item There is {\bf no} natural transformation $\rkw:\ID\to\KK$
non-vanishing on Verma modules.
\end{enumerate}
\end{proposition}

\begin{proof}
We consider the defining sequence
$0\rightarrow\KK\stackrel{i}\rightarrow\theta\stackrel{\ADJ^s}
\rightarrow\ID$. It induces an exact sequence
$\HOM(\ID,\KK)\stackrel{i\circ}\hookrightarrow\HOM(\ID,\theta)
\stackrel{\circ\ADJ^s}\rightarrow\HOM(\ID,\ID)$.
We have $\HOM(\ID,\theta)\cong \Coinv$, more precisely, the
morphism space is generated by the adjunction morphism $\ADJ_s$
and the center $\Coinv$ of the category $\cO_0$ (see \cite[Theorem~4.9]{Ba}).
If now $\phi\in\HOM(\ID,\KK)$ does not vanish on Verma modules, then,
up to a scalar, $i\circ\phi=\ADJ_s$, hence $\ADJ^s\circ i\circ\phi=
\ADJ^s\circ\ADJ_s\not=0$ (see \cite[Sections~2 and 3]{Be} or \cite[Lemma 2.2]{Afilt}).
This contradicts the exactness of the original exact sequence.
\end{proof}

\section{Proof of Theorem~\ref{tmain2}}\label{s3}

Theorem~\ref{tmain2} (\ref{tmain2.105}) follows
immediately from \cite[Section 4]{MS} and the definition of $\wZZ$.

\begin{proof}[Proof of Theorem~\ref{tmain2} (\ref{tmain2.1}).]
Let $\cH$ be the category of Harish-Chandra bimodules with
generalized trivial central character from both sides (see
\cite{SHC}). By \cite[5.9]{BG}, the category $\Oo_0$ is equivalent to the
full subcategory of $\cH$ given by objects having trivial central
character from the right hand side. Let $\theta_s^r:\cH\rightarrow\cH$
denote the right translation through the $s$-wall. When considering
$\cO_0$ as a subcategory of $\cH$, the functor $\GG$ is defined  as
the kernel of the adjunction morphism $\theta_s^r
\overset{\mathrm{adj}}{\longrightarrow}\ID$ (see \cite[Proposition~3.2]{J}).
Using the Snake Lemma we obtain that $\mathcal{R}^1\GG$ is
isomorphic to the cokernel of $\theta_s^r
\overset{\mathrm{adj}}{\longrightarrow}\ID$.
Note that $\cR^1G(M)$ is locally $\mg^s$-finite (\cite[Corollary
5.9]{AS}). Since the top of $\theta_s^r M$ is $s$-free, we obtain
that it is maximal with this property. Hence
$\cR^1G\cong\ZZ$ and, in particular, $\cR^1G\cong\ID$ on $\cO_0^s$.
\end{proof}

\begin{remark}{\rm
Theorem~\ref{tmain2}\eqref{tmain2.1} has independently been
proved in \cite[Proposition~20]{Kh} by completely different arguments.
\erem}
\end{remark}

\begin{proof}[Proof of Theorem~\ref{tmain2}\eqref{tmain2.2}.]
Recall from above that the functor $\ZZ$ is isomorphic to the cokernel
of the $\theta_s^r \overset{\mathrm{adj}}{\longrightarrow}\ID$. Let
$M\in\Oo_0$ and $P_2\overset{h}{\to} P_1\overset{f}{\to} P_0\tto M$ be
the first three steps of a projective resolution of $M$. Consider
the following commutative diagram:
\begin{displaymath}
\xymatrix{
\GG P_2\ar[rr] \ar@{^{(}->}[d]&& \GG P_1\ar[rr]\ar@{^{(}->}[d]&&
\GG P_0\ar[rr]\ar@{^{(}->}[d]&& \GG M\ar@{^{(}->}[d] \\
\theta_s^r P_2\ar[rr]\ar[d]^{\mathrm{adj}} &&
\theta_s^r P_1\ar[rr]\ar[d]^{\mathrm{adj}} &&
\theta_s^r P_0\ar@{->>}[rr]\ar[d]^{\mathrm{adj}} &&
\theta_s^r M\ar[d]^{\mathrm{adj}} \\
P_2\ar[rr]^{h}\ar@{->>}[d]^{p_2} &&
P_1\ar[rr]^{f}\ar@{->>}[d]^{p_1} &&
P_0\ar@{->>}[rr]\ar@{->>}[d]^{p_0} &&  M\\
\ZZ P_2  \ar[rr]^{\overline{h}} &&\ZZ P_1  \ar[rr]^{\overline{f}}
&&\ZZ P_0 &
}.
\end{displaymath}
The Snake Lemma gives a natural surjection $\GG M\surj
\ZZ (P_1/\Ker f)$. We claim that this even induces a natural
surjection $\GG M\surj\Ker\overline{f}/{\Im\overline{h}}$. Indeed,
if $x\in \ZZ P_1$ such that $\overline{f}(x)=0$ and
$x\not\in\Im\overline{h}$, we can choose $y\in P_2$ such that
$p_2(y)=x$. If  $f(y)=0$ then $y=h(z)$ for some $z\in P_3$; hence
$x=p_2\circ h(z)= \overline{h}\circ p_3(z)$, which is a contradiction.
Therefore, $f(y)\neq 0$ and $\ZZ (P_1/\Ker f)$ surjects onto
$\Ker\overline{f}/{\Im\overline{h}}$ providing a surjection $\Phi:\GG
\surj \mathcal{L}_1 \ZZ $. We have to show that $\Phi$ induces an
isomorphism $\FF\cong\cL_1\ZZ$.

\begin{claim}\label{l3.4}
\begin{displaymath}
\mathcal{L}_1\ZZ\Delta(x)\cong
\begin{cases}
\Delta(sx)/\Delta(x), &\text{if $x>sx$},\\
0, &\text{if $x<sx$}.
\end{cases}
\end{displaymath}
In particular, $\Phi$ induces an isomorphism
$\FF\cong\mathcal{L}_1\ZZ$ on Verma modules.
\end{claim}

\begin{proof}
We prove the claim by induction on $l(x)$. It is certainly true for
$x=e$. Assume it to be true for $x$ and let $t$ be a simple reflection
such that $xt>x$. The short exact sequence $\Delta(x)\hookrightarrow
\theta_t \Delta(x)\tto\Delta(xt)$ induces an exact sequence
\begin{equation}
\label{eq2}
\mathcal{L}_1\ZZ\Delta(x)\inj
\mathcal{L}_1\ZZ\theta_t\Delta(x)\to
\mathcal{L}_1\ZZ\Delta(xt)\to
\ZZ\Delta(x)\to\ZZ\theta_t \Delta(x)\tto
\ZZ\Delta(xt).
\end{equation}
{\it If $x>sx$} then $l(sxt)\leq l(sx)+1=l(x)<l(xt)$.
Since $x>sx$ and $sxt> xt$, we have
$\ZZ\Delta(x)=\ZZ\Delta(xt)=\ZZ\theta_t
\Delta(x)=0$. By induction
hypothesis, \eqref{eq2} reduces to
\begin{displaymath}
\Delta(sx)/\Delta(x)\inj \theta_t(\Delta(sx)/\Delta(x))\tto
\mathcal{L}_1\ZZ\Delta(xt),
\end{displaymath}
implying $\mathcal{L}_1\ZZ\Delta(xt)\cong \Delta(sxt)/\Delta(xt)$.

{\it If $sx>x$ and  $sxt<xt$} then $xt>x$ implies $sxt=x$.
Hence $\ZZ\Delta(xt)=\ZZ\theta_t\Delta(x)=\ZZ\theta_t\Delta(x)=0$,
and $\mathcal{L}_1\ZZ\theta_t\Delta(x)\cong\theta_t
\mathcal{L}_1\ZZ\Delta(x)=0$ by induction hypothesis. We get
\begin{displaymath}
\mathcal{L}_1\ZZ\Delta(xt)\cong \ZZ\Delta(x)\cong
\Delta(x)/\Delta(sx)=\Delta(sxt)/\Delta(xt).
\end{displaymath}

{\it If $sx>x$ and $sxt>xt$} then we have
$(\mathcal{L}_1\ZZ)\theta_t\Delta(x)\cong\theta_t
(\mathcal{L}_1\ZZ)\Delta(x)=0$ by induction hypothesis, and
the last terms of \eqref{eq2} form the exact sequence
\begin{displaymath}
\Delta(x)/\Delta(sx)\hookrightarrow
\theta_t\Delta(xt)/\Delta(sxt)\tto \Delta(xt)/\Delta(sxt).
\end{displaymath}
This implies that $\mathcal{L}_1\ZZ\Delta(xt)=0$ and the claim follows.
\end{proof}

\begin{claim}\label{l3.5}
$\Phi$ induces an isomorphism $\FF\cong\cL_1\ZZ$ on modules having
a Verma flag.
\end{claim}

\begin{proof}
Let ${\bf S}$ be a short exact sequence of modules having a
Verma flag; then we have a commutative diagram
${\bf S}\stackrel{\rg_{\bf S}}\inj
G({\bf S})\surj\FF({\bf S})\rightarrow\cL_1Z({\bf S})$,
where the composition of the last two maps is $\Phi$. Since $\rg$
is an injection, $\FF({\bf S})$ is left-exact by the Snake Lemma.
The sequence $\cL_2Z({\bf S})$ is identical zero, because
$\mathcal{L}_2\ZZ\cong\ZZ'$ by \cite[Theorem~4.3]{EW}. Therefore,
$\cL_1Z({\bf S})$ is left-exact. The Five-Lemma implies the claim.
\end{proof}

\begin{claim}\label{l3.6}
$\Phi$ induces an isomorphism $\FF\cong\cL_1\ZZ$ on modules having
a dual Verma flag.
\end{claim}

\begin{proof}
Let ${\bf S}$ be a short exact sequence of modules having a dual
Verma flag; then $G({\bf S})$ is exact (\cite[Theorem 2.2]{AS})
and hence $\FF({\bf S})$ is right exact. On the other hand
$\mathcal{L}_1\ZZ({\bf S})$ is right exact as well, since
$\ZZ M=0$ for any module having a dual Verma flag.
The Five-Lemma completes the proof.
\end{proof}

Let $M\in\Oo_0$. By Wakamatsu's Lemma (\cite[Lemma~1.2]{Wa}) there
exists a short exact sequence ${\bf S}: Y\hookrightarrow X\tto M$,
for a certain $X$ having a Verma flag and some $Y$ with a dual Verma
flag. Since $\mathcal{R}^1 G(Y)=0$ (\cite[Theorem~2.2]{AS}), the
sequence $G({\bf S})$ is exact, and hence $Q({\bf S})$ is right
exact. Since $\ZZ Y=0$, $\mathcal{L}_1\ZZ({\bf S})$ is right exact,
as well. The Five-Lemma implies that $\Phi$ induces an isomorphism
$\FF M\cong\cL_1Z M$. We immediately get $\FF\cong\FF'$, since
$\mathcal{L}_1\ZZ\cong (\mathcal{L}_1\ZZ)'$ by
\cite[Theorem~4.3]{EW}. Theorem~\ref{tmain2}\eqref{tmain2.2} follows.
\end{proof}

\begin{proof}[Proof of Theorem~\ref{tmain2} (\ref{tmain2.3}).]
Recall the isomorphism $\mathcal{R}^1\GG\cong \ZZ$ from the first
part. By \cite[Theorem~2.2 and Theorem~4.1]{AS}, we have $\mathcal{R}^i\GG=0$ for all $i>1$.
Since $G(\rd\Delta(e))$ is acyclic for $G$
(\cite[Theorems 2.2 and 2.3]{AS}), we have the Grothendieck spectral
sequence $\mathcal{R}^p \GG(\mathcal{R}^q \GG
(X))\Rightarrow\mathcal{R}^{p+q} \GG^2 (X)$.
We immediately get $\cR^1\GG^2\cong\cZ G$ and $\cR^2\GG^2\cong
\ZZ^2\cong\ZZ$ and $\cR^i \GG^2=0$ for $i>2$. This
proves the first part of Theorem~\ref{tmain2}\eqref{tmain2.3}.

The second part is proved by analogous arguments provided that we know
that $\KK(I)$ is $\KK$-acyclic for any injective object $I$. This is
equivalent to the statement that the head of $\KK(I)$ contains no composition
factor $L(w)$ with $ws>w$. There is a short exact sequence
$X\hookrightarrow Y\tto I$, where $X$ has a dual Verma flag and
$Y$ is the projective-injective cover of $I$.  Using that $\KK$ is exact on
sequences of modules having a dual Verma flag, we get a surjection $\KK(Y)\tto
\KK(I)$. In particular, it follows that the head of $\KK(I)$ is
embedded into the head of $\KK(Y)\in\mathrm{Add}(P(w_0))$.
The latter contains only copies of $L(w_0)$. This completes the
proof.
\end{proof}

\section{Proof of Theorem~\ref{tmain3}}\label{s5}

We start by verifying the indicated relations. By duality, it is
enough to prove every second statement.

{\it The isomorphism $\TT\GG\TT\cong\TT$}:
Evaluating the exact sequence of functors
\begin{equation}\label{eq3}
0\to \TT\GG{\hookrightarrow}\ID\tto\ZZ\to 0,
\end{equation}
from Proposition~\ref{l2.102}\eqref{1} at $\TT$ gives rise
to the exact sequence $0\to \TT\GG\TT{\hookrightarrow}
\TT\tto\ZZ\TT\to 0$. Further $\ZZ\TT=0$, as the head of any $\TT(M)$
is $s$-free by \cite[Corollary 5.2]{AS}, hence we obtain
$\TT\GG\TT\cong\TT$.

{\it The isomorphism $\GG^3\cong\GG^2$} is proved in \cite[Lemma~3.6]{J}.

{\it The isomorphism $\TT^2\GG\cong\TT^2$}: Applying $\TT$ to
\eqref{eq3} gives the exact sequence
\begin{equation}
\label{eq4}
(\mathcal{L}_1\TT)\ZZ\to
\TT^2\GG{\rightarrow}\TT\tto
\TT\ZZ\to 0.
\end{equation}
Theorem~\ref{tmain2} gives $\mathcal{L}_1\TT\cong \ZZ'$,
in particular, $\TT (\mathcal{L}_1\TT)\ZZ=0$ (\cite[Corollaries 5.8
and 5.9]{AS}). Moreover $\TT\ZZ=0$. This means that we can apply
$\TT$ to \eqref{eq4} once  more to obtain an isomorphism
$\TT^3\GG\cong\TT^2$. Since $\TT^3\cong \TT^2$ we finally get
$\TT^2\GG\cong \TT^2$.

{\it The isomorphism $\TT\GG^2\cong \GG\TT^2$}:
Evaluating the adjunction morphism $\ADJ_{\TT}:
\TT\GG\hookrightarrow\ID$ at
$\GG\TT^2$ we get
$\TT\GG\GG\TT^2\cong \TT\GG^2\hookrightarrow\GG\TT^2$.
Evaluating $\ID\tto \GG\TT$ at $\TT\GG^2$ we obtain
$\TT\GG^2\tto \GG\TT\TT\GG^2\cong \GG\TT^2$ and hence
$\TT\GG^2\cong \GG\TT^2$.

To complete the proof it is now enough to show that all the functors
from $\mathcal{S}$ are not isomorphic (Green's relation are easily
checked by direct calculations). An easy direct calculation gives the
following images under our functors:
\begin{displaymath}
  \begin{array}{c|c|c|c|c|c|c|c}
\ID&\GG&\TT&\GG^2&\TT^2&\TT\GG&\GG\TT&\GG\TT^2\\
\hline\hline
\Delta(s)&\Delta(e)&\TT\Delta(s)&
\Delta(e)&\TT\Delta(s)&\Delta(s)&\Delta(s)&\Delta(s)
\\
\hline
\Delta(e)&\Delta(e)&\Delta(s)&
\Delta(e)&\TT\Delta(s)&\Delta(s)&\Delta(e)&\Delta(s)\\
\hline
\TT\Delta(s)&\Delta(s)&\TT
\Delta(s)&\Delta(e)&\TT\Delta(s)&
\TT\Delta(s)&\Delta(s)&\Delta(s)
\end{array}
\end{displaymath}
The claim follows.

\section{Proof of Theorem~\ref{tmain3new}}\label{s501}

By duality it is enough to prove every second relation.

{\it The isomorphism $\CC\KK\CC\cong\CC$:} The proof is analogous to
that of $\TT\GG\TT\cong\TT$ in Section~\ref{s5}.

{\it The isomorphism $\CC^3\KK\cong\CC^2$}: Applying  $\CC$ to the
short exact sequence $\CC\KK\hookrightarrow \ID\tto\hat{Z}$ produces a
short exact sequence $X\hookrightarrow \CC^2\KK\tto\CC$, where
$\CC X=0$. Applying $\CC$ once more we obtain the desired isomorphism.

{\it The isomorphism $\CC^2\KK^2\CC\cong\CC^2\KK$}: Applying  $\KK$
to the short exact sequence $\hat{Z}'\hookrightarrow \ID\tto\KK\CC$
produces a short exact sequence $\KK\hookrightarrow \KK^2\CC\tto X$,
where  $\KK X=\CC X=0$. Applying  now $\CC$ gives rise to
$Y\hookrightarrow \CC\KK\tto \CC\KK^2\CC$, where  $\KK Y=\CC
Y=0$. Applying $\CC$ once more gives the isomorphism.

{\it The isomorphism $\KK\CC^2\KK^2\cong\CC\KK^2$}: Evaluating the
short exact sequence $\hat{Z}'\hookrightarrow \ID\tto\KK\CC$ at
$\CC\KK^2$ we obtain the short exact sequence
$\hat{Z}'\CC\KK^2\hookrightarrow \CC\KK^2\tto\KK\CC^2\KK^2$. The
statement follows if we show that $\hat{Z}'\CC\KK^2=0$. The injection
$\CC\KK\hookrightarrow \ID$ gives an injection
$\CC\KK^2\hookrightarrow \KK$. On the other hand, $\hat{Z}'\KK=0$
since, by the definition of $\KK$, any composition factor in the socle
of $\KK M$ is not annihilated by $\theta$. As $\CC\KK^2\hookrightarrow
\KK$ we get that $\hat{Z}'\CC\KK^2=0$ as well.

It is easy to see that, using the relations we have just proved,
any product of $\CC$ and $\KK$ can be reduced to one of the elements
of $\mathcal{\hat{S}}$.

Assume now that $s$ does not correspond to an $\mathfrak{sl}_2$-direct
summand of $\mathfrak{g}$. We do a case-by-case analysis to show that
all functors in $\mathcal{\hat{S}}$ are different. We start with the
following general observation.

\begin{lemma}\label{cln5.1}
Assume that $\XX:\Oo_0\to\Oo_0$ is left exact, $\XX(P(w_0))\cong
P(w_0)$, and there is a natural transformation $\varphi:\ID\to \XX$
on the category of projective-injective modules in $\Oo_0$,
such that $\varphi_{P(w_0)}$ is an isomorphism. Then $\XX$ fixes the
isoclasses of projectives.
\end{lemma}

\begin{proof}
Let $P$ be projective. Consider an exact sequence  $P\hookrightarrow
I_0\to I_1$, where $I_0$ and $I_1$ are projective-injective. Then the
square on the right hand side in the following diagram with exact rows
commutes
\begin{displaymath}
\xymatrix{
 0\ar[rr] && P \ar[rr]^f\ar@{-->}[d]^{h}&& I_0 \ar[rr]^g
\ar[d]^{\varphi_{I_0}}&& I_1\ar[d]^{\varphi_{I_1}} \\
 0\ar[rr] && \XX P \ar[rr]^{\XX(f)}&& \XX I_0 \ar[rr]^{\XX(g)}&&
\XX I_1 \\
}
\end{displaymath}
and hence we obtain the induced map $h$, which is an isomorphism by
the Five Lemma.
\end{proof}

{\it All $\KK^i$ are different.}
We fix a simple reflection $t$ such that $st\neq ts$. By a direct
calculation one obtains that $\KK^i P(t)$, $i>0$, is not projective,
in particular, $\KK^i$ does not preserve projectives in $\Oo_0$. Now
any isomorphism $\varphi:\KK^i\to \KK^j$, $i<j$, induces a natural
transformation $\ID\to \KK^{j-i}$ on the category $\KK^i(\Oo_0)$,
which contains the subcategory of projective-injective modules in
$\Oo_0$.  It follows from Lemma~\ref{cln5.1} that $\KK^{j-i}$
preserves the category of projective modules in $\Oo_0$, a
contradiction.

{\it All $\CC^i$ are different} by dual arguments.

We consider now $\mathcal{\hat{S}}$ as a $\Z$-graded monoid with
$\mathrm{deg}(\CC)=1$ and $\mathrm{deg}(\KK)=-1$. This is possible as
the defining relations are homogeneous with respect to this
grading. It follows from the relations that for any $X\in
\mathcal{\hat{S}}$ and for all $i$ large enough we have
$\CC^i X\cong\CC^j$ for some $C^j$. Since we have already shown
that all $\CC^j$ are different, it follows that the elements of
$\mathcal{\hat{S}}$ having different degree are not isomorphic. In
particular, changing the exponent $i$ in the expression for
$X\in\mathcal{\hat{S}}$ gives a non-isomorphic functor.  The rest
will be checked case-by-case.

{\it $\KK^{i}$ is not isomorphic to $\CC\KK^{i+1}$ for $i>0$}: We have
$\CC\KK^{i+1}\Delta(e)\cong \Delta(s)$ and $\KK^{i} \Delta(e)\cong
\Delta(e)$ for all $i$.

{\it $\KK^{i}$ is not isomorphic to $\CC^2\KK^{i+2}$ for $i>0$}: We
have $\KK^{i+2} \Delta(e)\cong\Delta(e)\not\cong\CC\Delta(s)\cong
\CC^2\KK^{i+2} \Delta(e)$.

{\it $\KK$ is not isomorphic to $\KK^2\CC$}, since
$\KK\rd\Delta(e)\not\cong\KK^2\rd\Delta(e)\cong\KK^2\CC\rd\Delta(e)$.

We proved that $\KK^i$ (where $i>0$) is not isomorphic to any other
functor in the list. By duality, the same holds for $\CC^i$.

{\it $\KK\CC$ is not isomorphic to $\CC\KK$}: Assume, they are
isomorphic, then $\CC\cong\CC\KK\CC\cong\CC\CC\KK\cong\CC^2\KK$ which
we have proved to be wrong.

{\it $\KK\CC^i$ is not isomorphic to $\KK^2\CC^{i+1}$ for $i>0$}:
We have $\KK\CC^i\rd\Delta(e)\cong\KK\rd\Delta(e)\not\cong\KK^2\rd
\Delta(e)\cong\KK^2\CC^{i+1}\rd\Delta(e)$.

{\it $\KK\CC^2$ is not isomorphic to $\CC^2\KK$}: We have
$\KK\CC^2\rd\Delta(e)\cong\KK\rd\Delta(e)\cong\rd\Delta(s)$ and
$\CC^2\KK\rd\Delta(e)\cong\CC^2\rd\Delta(s)\cong\CC\rd\Delta(e)
\cong\rd\Delta(e)$.

{\it $\KK\CC$ is not isomorphic to $\KK\CC^2\KK$}: Assume, they are
isomorphic. Then $\KK\cong\KK\CC\KK\cong\KK\CC^2\KK^2\cong\CC\KK^2$,
which we know is wrong.

Hence the functors $\KK\CC^i$, $i>0$, differ from all the others in
the list. Duality gives the same property for $\CC\KK^i$.

{\it $\KK^2\CC^2$ is not isomorphic to $\CC^2\KK^2$ and $\KK^2\CC$ is
not isomorphic to $\CC^2\KK^3$}: By definition the socle of
$\KK^2\CC^2 M$ contains only composition factors which are not
annihilated by $\theta$ (for any $M\in\cO_0$). On the other hand
$\CC^2\KK^2\Delta(e)\cong\CC^2\Delta(e)\cong\CC\Delta(s)$ is an
extension of $\Delta(s)$ with $\Delta(e)/\Delta(s)$. In particular,
the socle is $\mg^s$-finite. The same argumentation applies to the
second pair.

{\it $\KK^2\CC^2$ is not isomorphic to $\KK\CC^2\KK$}: Assume, they
are isomorphic then $\KK^2\CC\cong\KK^2\CC^2\KK\cong\KK\CC^2\KK^2
\cong\CC\KK^2$. We have already proved that this is not possible.

Hence $\KK^2\CC^i$, $i>0$, (and dually $\CC^2\KK^i$) differs from all
other functors from the list. And therefore, any two functors from
the list are not isomorphic.

The statements concerning Green's relations and idempotents
are obtained by a direct calculation.

\section{Proof of Theorem~\ref{tmain1}}\label{s4}

It will be enough to prove roughly
half of the statements.  The other half will follow by duality.

\begin{lemma}\label{pr4.1}
All maps indicated in the diagram as inclusions are injective;
and all projections are surjective.
\end{lemma}

\begin{proof}
By duality, it is enough to prove the statement for inclusions.
The injectivity of $\rz'$, $\ri'$, $\ri'_T$, $\rz'_T$, $\rz'_{T^2}$
is given by definition. For the maps $\GG(\rt)$ and  $\GG(\rg)$ the
statement follows from the left exactness of $\GG$ and the fact that
$\GG$ is zero on locally $\mg^s$-finite modules. The map
$\ZZ'(\ri_{\TT})$ is injective because of the left exactness of $\ZZ'$
and the injectivity of $\ri_{\TT}$. The injectivity of $\raa'$ follows
from \cite[Proposition~5.6]{AS}, since $\raa'$ is up to a non-zero
scalar the adjunction morphism $\ADJ_{\TT}:\TT\GG\to\ID$.

Let us now prove the statement for $\ZZ\GG(\rg)$. By definition of
$\FF$ we have the following exact sequence of functors:
$\GG\hookrightarrow\GG^2\tto \FF\GG$. It gives rise to the exact
sequence
\begin{displaymath}
0\cong \mathcal{L}_1\ZZ(\FF\GG)\to \ZZ\GG
\overset{\ZZ\GG(\rg)}{\longrightarrow}
\ZZ G^2\overset{\GG(\rp_{\GG})}{\tto}\ZZ\FF\GG\cong \FF\GG.
\end{displaymath}
This implies that $\ZZ\GG(\rg)$ is injective.

\begin{claim}\label{l4.2}
$\TT^2(\rg):\TT^2\to\TT^2\GG$ is an isomorphism. In particular
$\rk'$ is well-defined and injective.
\end{claim}

\begin{proof}
Let $K$ and $K'$ be defined by the following exact sequence of
functors:
\begin{displaymath}
\xymatrix{
K\ar@{^{(}->}[rr] && \ID\ar[rr]^{\rg}\ar@{->>}[rd]^{q} &&
\GG\ar@{->>}[rr] && K'\\
 && &\mathrm{im}(\rg)\ar@{^{(}->}[ru]^{j}& &&
},
\end{displaymath}
Since $\TT^2 K=0$ we get an isomorphism
$\TT^2(q):\TT^2\to \TT^2(\mathrm{im}(\rg))$ where $\mathrm{im}(\rg)$
denotes the image of of $\rg$. Applying $\TT$ to the second short
exact part  gives a short exact sequence $\tilde{K}\hookrightarrow
\TT(\mathrm{im}(\rg))\tto \TT\GG$ for some $\tilde{K}$ such that
$\tilde{K}(M)$ is locally $\mg^s$-finite for all $M\in\cO_0$. Applying
$\TT$ once more gives an isomorphism $\TT^2(j):\TT^2(\mathrm{im}(\rg))
\to \TT^2\GG$ since $\TT\tilde{K}=0$. Composing $\TT^2(j)\circ
\TT^2(q)=\TT^2(\rg)$ implies the first statement. The injectivity of
$\rk'$ follows from the injectivity of $\ri_{\TT\GG}$.
\end{proof}

\begin{claim}\label{l4.201}
There exists a unique isomorphism $\rh:\TT\GG^2\to\GG\TT^2$ such that
\begin{displaymath}
\rg\circ\rt=\GG(\rt\circ\rt_{\TT})\circ\rh\circ\TT(\rg_{\GG}\circ\rg).
\end{displaymath}
\end{claim}

\begin{proof}
We start proving uniqueness. If
$\rh$ and $\tilde{\rh}$ are two such morphisms, then
$\rh-\tilde{\rh}$ induces a morphism from $\ZZ'\TT$ to $\GG$ since
$\ZZ'\TT=\ker (\rg\circ\rt)$ (this will be proved later in this
section). However, $\Hom(\ZZ'\TT,\GG)=0$ as the socle of $\GG M$ is
$s$-free and $\ZZ'\TT M$ is $\mg^s$-finite for any $M\in\Oo_0$.

It is left to prove the existence. Note that
$\TT\GG^2\cong\GG\TT^2$ by Theorem~\ref{tmain3}.
For any $h\in \End(\TT\GG^2,\GG\TT^2)$ the natural transformation
$\varphi(h)=\GG(\rt\circ\rt_{\TT})\circ h\circ\TT(\rg_{\GG}\circ\rg)$
belongs to $\Hom(\TT,\GG)$ and, comparing the action on the
projective-injective module $P(w_0)\in\Oo_0$ we see that
$\varphi$ is injective, hence an isomorphism (by the independent
Theorem~\ref{tmain4}). The claim follows.
\end{proof}

We proceed with the map $\FF'\TT(\rg)$. Let $M\in\Oo_0$ and consider
the map $\rg_M:M\to\GG M$. The map $T(\rg_M)$ fits into the exact
sequence $\FF' M\to \TT M\to \TT\GG M$. To calculate $\FF'\TT(\rg)$
we consider the following commutative diagram:
\begin{displaymath}
\xymatrix{
\FF'\FF' M=0\ar[rr]\ar@{^{(}->}[d] &&
\FF' \TT M\ar[rr]^{\FF'\TT(\rg_M)}\ar@{^{(}->}[d]
&& \FF'\TT\GG(M)\ar@{^{(}->}[d]\\
\TT\FF' M=0\ar[rr]\ar[d] && \TT^2 M\ar[rr]^{T^2(\rg_M)}\ar[d] &&
\TT^2\GG M\ar[d]\\
\TT\GG\FF' M=0\ar[rr] && \TT\GG\TT M\ar[rr]^{\TT\GG\TT(\rg_M)}&&
\TT\GG\TT\GG M \\
},
\end{displaymath}
where the first row is the kernel sequence and hence is exact. It
follows that $\FF'\TT(\rg)$ is injective. The injectivity of
$\FF(\rg\circ\rt)$ is proved by analogous arguments. This completes
the proof of Lemma~\ref{pr4.1}.
\end{proof}

\begin{lemma}\label{pr4.5}
All configurations containing only solid arrows
commute.
\end{lemma}

\begin{proof}
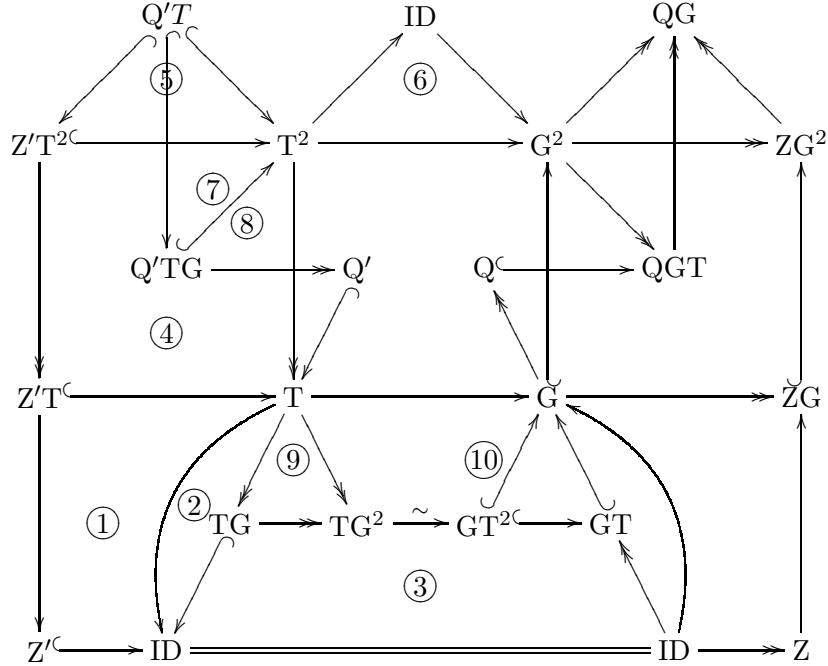
\begin{figure}[htbp]
  \centering
$\xymatrix@!=0.001pc{
 &&\FF'T\ar@{^{(}->}[lldd]
\ar@{^{(}->}[rrdd]
\ar@{^{(}->}[dddd]
&& && \ID
\ar[ddrr]
&& && \FF\GG && \\
&& *+[o][F-]{5}&& &&*+[o][F-]{6}\\
\ZZ'\TT^2 \ar@{->>}[dddd]
\ar@{^{(}->}[rrrr]
&& && \TT^2 \ar[rrrr]
\ar[rruu]\ar@{->>}[dddd]
&& && \GG^2 \ar@{->>}[rrrr]
\ar@{->>}[rruu]\ar@{->>}[rrdd]
&& && \ZZ\GG^2\ar@{->>}[lluu] \\
\\
 && \FF'\TT\GG \ar@{^{(}->}[rruu]^*+[o][F-]{7}_*+[o][F-]{8}
\ar@{->>}[rrr]
& && \FF'\ar@{^{(}->}[ddl]
 && \FF \ar@{^{(}->}[rrr] & &&
\FF\GG\TT
\ar@{->>}[uuuu] && \\
&&*+[o][F-]{4}\\
\ZZ'\TT\ar@{^{(}->}[rrrr]
\ar[dddd]
 && && \TT \ar[rrrr]
\ar@/_2pc/[ddddll]^*+[o][F-]{2}\ar@{->>}[ldd]\ar@{->>}[rdd]
&& && \GG \ar@{->>}[rrrr]
\ar@{_{(}->}[uuuu]\ar@{->>}[luu]
&& && \ZZ\GG\ar@{^{(}->}[uuuu] \\
&& && *+[o][F-]{9}& && *+[o][F-]{10}\\
 & *+[o][F-]{1}&&  \TT\GG\ar@{^{(}->}[ldd]
\ar@{->>}[rr]  & &
\TT\GG^2\ar[rr]^{\sim}
& &\GG\TT^2\ar@{^{(}->}[rr]\ar@{^{(}->}[ruu]& &
\GG\TT\ar@{_{(}->}[luu]\\
&& && &&*+[o][F-]{3}\\
\ZZ'\ar@{^{(}->}[rr]
&&
\ID
\ar@{=}[rrrrrrrr]
& && && && &\ID
\ar@{->>}[rr]\ar@{->>}[luu]
\ar@/_2pc/[uuuull]
&& \ZZ \ar[uuuu]\\
}$
\caption{Schematic picture of the diagram from Theorem~\ref{tmain1}}
  \label{P1}
\end{figure}
We use the notations from Figure~\ref{P1}.
The squares $\xymatrix{*+[o][F-]{2}}$,
$\xymatrix{*+[o][F-]{6}}$,
$\xymatrix{*+[o][F-]{9}}$, and
$\xymatrix{*+[o][F-]{10}}$ commute by definition. The commutativity of
$\xymatrix{*+[o][F-]{3}}$ follows from the commutativity
of $\xymatrix{*+[o][F-]{2}}$, $\xymatrix{*+[o][F-]{9}}$, and
$\xymatrix{*+[o][F-]{10}}$. The squares $\xymatrix{*+[o][F-]{1}}$, and
$\xymatrix{*+[o][F-]{4}}$ commute since $\rz'$ is a natural
transformation and $\ZZ'$ and $\ZZ'\TT$ are functors
(note that $\rt_{\TT}=\TT(\rt)$). The commutativity of
$\xymatrix{*+[o][F-]{5}}$ reads $\ri_T=\rz'_{\TT^2}\circ \ZZ'(\ri_T)$,
which is true as $\ZZ'=\ID$ on $\mg^s$-finite modules. The
commutativity of $\xymatrix{*+[o][F-]{7}}$ reads
$\ri_T=\rk'\circ \FF'\TT(\rg)$, which is equivalent to
$\TT^2(\rg)\circ \ri_{\TT}=\ri_{\TT\GG}\circ \FF'\TT(\rg)$, the latter
being true as $\ri$ is a natural transformation. Commutativity of
$\xymatrix{*+[o][F-]{8}}$ means $\ri\circ
\FF'(\raa')=\rt_{T}\circ\rk'$, which is equivalent to
$\ri\circ \FF'(\raa')=\rt_{T}\circ (\TT^2(\rg))^{-1}\circ\ri_{\TT\GG}$.
Since $\ri$ is a natural transformation we have
$\ri\circ \FF'(\raa')=\TT(\raa')\circ \ri_{\TT\GG}$ and our equality
reduces to $\TT(\raa')\circ \ri_{\TT\GG}=
\rt_{T}\circ (\TT^2(\rg))^{-1}\circ\ri_{\TT\GG}$. To prove the latter
it is enough to show that $\TT(\raa')=\rt_{T}\circ (\TT^2(\rg))^{-1}$,
which follows from $\rt_{T}=\TT(\rt)$ and the definition of $\raa'$.
The remaining configurations commute by duality.
\end{proof}

To complete the proof of Theorem~\ref{tmain1} it is left to prove the
exactness of the indicated sequences. By duality, it is sufficient to
prove the exactness of the sequences {\bf 1} to {\bf 10}.
The sequences {\bf 8} and {\bf 3} are exact by the
definitions of $\raa$ and $\FF$ respectively.
The exactness of {\bf 4} follows from \cite[Proposition~5.6]{AS}.
The exactness of {\bf 7} follows from $\TT(\rt)=\rt_{\TT}$ and
the exactness of the sequence, dual to {\bf 3}.
Applying the left exact functor $\ZZ'$ to the short exact sequence
{\bf 7} and using $\ZZ'\FF'=\FF'$ shows that {\bf 5}
is exact. The exactness of {\bf 6} follows by comparison of
characters from the facts that $\FF'\TT(\rg)$ is an inclusion and
$\FF'(\raa')$ is a surjection. The exactness of {\bf 10}
follows by evaluating the exact sequence {\bf 8} at modules of
the form $\GG M$.

Let us now show that {\bf 2} is exact. The cokernel
$\mathrm{Coker}$ of $\rg\circ\rt:\TT\to\GG$ is $\mg^s$-finite
since already the cokernel of $\rg$ is $\mg^s$-finite, see \cite[Lemma~3.10]{J}.
Further, for any $M\in\Oo_0$ we have that $\FF(M)$
is the maximal $\mg^s$-finite quotient of $\GG M$ since the head of
$\TT M$ is $s$-free. This implies the exactness of the sequence
{\bf 2} and also of {\bf 9} at the term $\GG$. By uniqueness of the
canonical maps the exactness in $\TT$ follows by duality.
Exactness of {\bf 1} follows by analogous arguments.

\section{Proof of Theorem~\ref{tmain4}}\label{s6}

We abbreviate $\Hom(X,Y)=\ho_{X,Y}$ for
$X,Y\in\mathcal{S}$. By duality we have vector space isomorphisms
$\ho_{X,Y}\cong \ho_{Y',X'}$.

\begin{proposition}\label{pr6.1}
$\End(X)\cong \Coinv$ as algebras for any
$X\in\mathcal{S}$.
\end{proposition}

\begin{proof}
For $X=\ID$ the statement is well-known and follows from
\cite[Endomorphismensatz and Struktursatz]{S}, since $\End(\ID)\cong \Coinv\cong\END_\mg(P(w_0))$.
Note that $\GG P(w_0)\cong \TT P(w_0)\cong P(w_0)$ (see
\cite[Proposition 5.3]{AS}); hence
$X P(w_0)\cong P(w_0)$ for all $X\in\mathcal{S}$. This means that
sending $\varphi\in \End(\ID)$ to $X(\varphi)$ defines an injective
algebra morphism from $\Coinv$ to $\End(X)$ for every
$X\in\mathcal{S}$, as already the map $\varphi_{P(w_0)}\mapsto
X(\varphi_{P(w_0)})$ is injective. We only have to check the
dimensions.

We claim that $\Phi:\End(\TT)\rightarrow\END_\mg(\TT P(w_0))$,
$\phi\mapsto\phi_{P(w_0)}$, is injective. Assume that
$\Phi(\phi)=0$. Let $P\in\cO_0$ be projective with injective hull
$i:P\hookrightarrow I$. The cokernel $Q$ has a Verma flag, hence
$0\rightarrow\TT P\stackrel{\TT i}{\inj}\TT I\surj\TT Q\rightarrow 0$
is exact (see \cite[Theorem 2.2]{AS}). Since $I$ is a direct sum of
copies of $P(w_0)$, we have $\phi_I=0$ and therefore $\phi_P=0$.
Since $\TT$ is right exact we get $\phi_M=0$ for any $M\in\cO_0$.
Hence $\Phi$ is injective and $\END(\TT)\cong \Coinv$. We get
$\END(\GG)\cong \Coinv$ by duality.

The adjointness from Proposition~\ref{l2.102} together with
Theorem~\ref{tmain3} imply
\begin{gather*}
\END(\TT^2)\cong\HOM(\ID,\GG^2\TT^2)\cong\HOM(\ID,\GG^2)\cong
\END(T)\cong \Coinv,\\ \End(\GG\TT)\cong \Hom(\TT\GG\TT,\TT)\cong
\End(\TT)\cong \Coinv
\end{gather*}
and also $\End(\GG\TT^2)\cong
\Hom(\TT\GG\TT^2,\TT^2)\cong \End(\TT^2)\cong \Coinv$. The
remaining parts follow by duality.
\end{proof}

\begin{claim}\label{homneq0}
$\ho_{X,Y}\not=0$ for any $X$, $Y\in\cS$.
\end{claim}

\begin{proof}
Since both $X$ and  $Y$ are isomorphic to the identity functor
when restricted to $\cA=\mathrm{Add}(P(w_0))$ (see Lemma~\ref{l2.1})
we can fix a natural transformation $\phi\in
\Hom(X|_{\cA},Y|_{\cA})\cong
\Coinv$ of maximal degree. For $M\in\cO_0$ indecomposable, $M\notin\cA$,
we set $\phi_M=0$. For $M\in\cO_0$ arbitrary we fix an isomorphism
$\alpha_M: M\cong M_1\oplus M_2$, such that $M_1$ is a maximal direct
summand belonging to $\cA$ and set $\phi_M:=X(\alpha_M^{-1})\circ
(\phi_{M_1}\oplus\phi_{M_2})\circ X(\alpha_M)$. We claim that this
defines an (obviously nontrivial) element $\phi\in\ho_{X,Y}$. Indeed,
let $M\cong M_1\oplus M_2$ and $N\cong N_1\oplus N_2$ and
$f\in\HOM_\mg(M,N)$ with decomposition $f=\sum_{i,j=1}^{2}f_{i,j}$
such that $f_{i,j}\in\HOM_\mg(M_i,N_j)$. Then
$\phi_N\circ X(f_{1,1})=Y(f_{1,1})
\circ\phi_M$ by definition of $\phi$. The definitions also immediately
imply $0=Y(f_{2,2})\circ\phi_M=\phi_N\circ X(f_{2,2})$. Moreover, we
also have $0=\phi_N\circ X(f_{1,2})$ and $0=Y(f_{2,1})\circ\phi_M$.
Indeed, if
$Y(f_{1,2})\circ\phi_M\not=0$ or $\phi_N\circ X(f_{1,2})\not=0$ then
either a direct summand of $Y(M_1)$ embeds into $Y(N_2)$ or $X(M_2)$
surjects onto a direct summand of $Y(N_1)$. Both contradict the
following statement: Assume $R\in\cS$ and $M\in\cO_0$ does not have
$P(w_0)$ as a direct summand then neither so does $R(M)$. Let first
$R\in\{\GG,\CC\}$. If $P(w_0)$ is a direct  summand of $R(M)$ then
$R'R M$ surjects onto $R' P(w_0)\cong P(w_0)$, hence $P(w_0)$ is a
direct summand of $ R'R M$. The inclusion $R'R\inj\ID$ from
Proposition~\ref{l2.102} implies that $P(w_0)$ is a submodule (hence a
direct summand) of $M$. Dual arguments apply to $R\in\{\TT,\KK\}$
and the claim follows.
\end{proof}

\begin{claim}\label{pr6.2}
The $\Coinv$-entries in the table of Theorem~\ref{tmain4} are correct.
\end{claim}

\begin{proof}
The statement is obtained by playing with the adjointness of $\TT$ and
$\GG$ using Proposition~\ref{pr6.1} and the identities from
Theorem~\ref{tmain3}. Let $X$, $Y\in\mathcal{S}$. We have isomorphisms
$\ho_{\TT^2,X}\cong \ho_{\TT^2\GG^2,X}\cong
\ho_{\GG^2,\GG^2X}\cong \ho_{\GG^2,\GG^2}\cong \Coinv$.  This gives the
spaces in question in the seventh row (and the sixth column by duality).
The isomorphisms $\ho_{\TT\GG,\ID}\cong \ho_{\GG,\GG}\cong \Coinv$ and
$\ho_{\TT\GG,\GG X}\cong\ho_{\TT^2\GG,X}\cong\ho_{\TT^2,X}\cong \Coinv$
imply the claim for the fifth row (and the fourth column by
duality). The spaces in question in the first, third and fourth rows
follow from $\ho_{\TT X,\GG Y}\cong\ho_{\TT^2 X, Y}\cong \Coinv$ and
$\ho_{\GG\TT,\GG}\cong\ho_{\TT\GG\TT,\ID}\cong\ho_{\TT,\ID}
\cong\ho_{\ID, \GG}, \ho_{\ID, \GG\TT\GG}\cong\ho_{\TT, \TT\GG}$.
From $\ho_{\GG\TT^2,\GG}\cong\ho_{\TT\GG\TT^2,\ID}\cong
\ho_{\TT^2,\ID}\cong \Coinv$ and $\ho_{\GG\TT^2,\GG\TT^2}\cong
\ho_{\TT\GG\TT^2,\TT^2}\cong\ho_{\TT^2,\TT^2}\cong \Coinv$ we get the
spaces in the last row. This completes the proof.
\end{proof}

To proceed we use the following general statement:

\begin{proposition}\label{p6.4}
Let $\mathfrak{A}$ be an abelian category with enough projectives. Let
$F$, $J$, $H$ be endofunctors on $\mathfrak{A}$. Assume that $F$
preserves surjections, and for any projective $P\in \mathfrak{A}$
there exists some $N\in \mathfrak{A}$ such that $F(P)\cong FH(N)$.
Then the restriction defines an injective map
$\Hom(F,J)\hookrightarrow \Hom(FH,JH)$.
\end{proposition}

\begin{proof}
It is enough to show that for any $\varphi\in \Hom(F,J)$ such that
$\varphi_{H}=0$ we have $\varphi=0$.
Let $M\in \mathfrak{A}$ with projective cover $f:P\tto M$. We choose
$N\in \mathfrak{A}$ such that $F(P)\cong FH(N)$. The first row of the
following commutative diagram is exact, since $F$ preserves
surjections.
\begin{displaymath}
\xymatrix{
FH(Q)\cong FP\ar@{->>}[rr]^{f}\ar[d]_{\varphi_{H(Q)}} &&
F(M)\ar[d]_{\varphi_{M}}\ar[r] & 0 \\
JH(Q) \ar@{->>}[rr] &&  G M
}.
\end{displaymath}
The surjectivity of $f$ and $\varphi_{H(Q)}=0$ imply $\varphi_{M}=0$.
\end{proof}

{\it The spaces with labeling different from $4$:} The indicated
equalities with labeling different from $1$ and $4$ follow directly
by duality. By \cite[Corollary~4.2]{AS}, the adjunction morphism
$\ADJ^{\TT}:\ID\tto \GG\TT(P)$
is an isomorphism on projectives. Hence, we may apply
Proposition~\ref{p6.4} to  $F=\ID$, $J=\TT$, and $H=\GG\TT$ to obtain
$\ho_{\ID,\TT}\hookrightarrow \ho_{\GG\TT,\TT\GG\TT}\cong
\ho_{\GG\TT,\TT}$. Further, the adjunction morphism
$\ADJ_{\TT}:\TT\GG\inj\ID$ is injective,  hence
$\ho_{\GG,\TT\GG}\inj\ho_{\GG,\ID}$ and
$\ho_{\GG\TT,\TT}\inj\ho_{\ID,\TT}$ by duality.

{\it The equality of the spaces denoted by $4$:} we have the following
isomorphisms
\begin{eqnarray}
\label{A}
  &\ho_{\GG\TT^2,\TT\GG}\cong\ho_{\TT\GG^2,\TT\GG}
\cong\ho_{\GG^2,\GG\TT\GG}\cong\ho_{\GG^2,\GG}&\\
\label{B}
&\ho_{\GG,\GG\TT^2}\cong\ho_{\TT\GG^2,\TT}\cong
\ho_{\GG^2,\GG\TT}\cong\ho_{\TT\GG,\TT^2}&\\
\label{C}
&\ho_{\GG^2,\GG\TT^2}\cong\ho_{\TT\GG^2,\TT^2}
\cong\ho_{\GG\TT^2,\TT^2}&\\
\label{D}
&\ho_{\GG,\GG\TT}\cong\ho_{\TT\GG,\TT}&.
\end{eqnarray}
Note that all the spaces labeled by $4$ occur in this list. The
inclusion $\TT\GG\inj\ID$ provides inclusions
$\GG\TT^2\cong\TT\GG^2\inj\GG$ and
$\TT\GG^2\cong\TT\GG^2\TT\inj\GG\TT$; hence
$\ho_{\GG^2,\GG\TT^2}\inj\ho_{\GG^2,\GG}$ and
$\ho_{\GG,\TT\GG^2}\inj\ho_{\GG,\GG\TT}$ (i.e. \eqref{C} is `included'
in \eqref{A} and \eqref{B} is `included' in \eqref{D}). Applying
Proposition~\ref{p6.4} with $F=\GG\TT^2$, $J=\TT$ and
$H=\TT$ ($F=\ID$, $J=\GG\TT^2$, $H=\GG$ respectively) we get
inclusions $\ho_{\GG\TT^2,\TT}\inj\ho_{\GG\TT^2,\TT^2}$ and
$\ho_{\ID,\GG\TT^2}\inj\ho_{\GG,\GG\TT^2\GG}\cong\ho_{\GG,\GG\TT^2}$
(i.e. \eqref{B} is `included' in \eqref{C} and \eqref{A} is `included'
in \eqref{B}). Hence, all the spaces from \eqref{A}--\eqref{D} have
the same dimension.

{\it The existence of the inclusions from {\bf A}}:
The inclusion $\TT\GG\inj\ID$ implies
$\ho_{\GG\TT,\TT\GG}\inj\ho_{\GG\TT,\ID}$.  Applying
Proposition~\ref{p6.4} to $F=\ID$, $J\in\{\TT,\TT\GG\}$, and
$H=\GG^2$, (this is possible since $\GG^2(P)\cong P$ for any
projective $P$) we get inclusions $\ho_{\ID,\TT}\hookrightarrow
\ho_{\GG^2,\TT\GG^2}$ and $\ho_{\ID,\TT\GG}\hookrightarrow
\ho_{\GG^2,\TT\GG^2}$. Finally, the inclusion $\GG\inj\GG^2$ gives
$\ho_{\GG^2,\GG}\inj\ho_{\GG^2,\GG^2}\cong \Coinv$.

{\it The existence of the inclusions from {\bf B}}:
Applying Proposition~\ref{p6.4} to $F=\ID$, $J=\TT^2$ and
$H\in\{\GG,\GG^2\}$, we obtain the inclusions
\begin{eqnarray}\label{eq12}
\ho_{\ID,\TT^2}\hookrightarrow \ho_{\GG,\TT^2},
&&
\ho_{\ID,\TT^2}\hookrightarrow \ho_{\GG^2,\TT^2}.
\end{eqnarray}
Finally, using again the adjunction $\TT\GG\inj\ID$ we get
$\ho_{\GG^2,\TT\GG}\inj\ho_{\GG^2,\ID}$.

{\it The existence of the inclusion {\bf C}}: We use the
following result (which generalizes without problems to arbitrary
parabolic subalgebras):

\begin{proposition}\label{l6.25}
There is a natural isomorphism $\End(\ZZ)\cong \END(\ID_{\cO_0^s})$ of rings.
\end{proposition}

\begin{proof}
Denote by $\PT$ the direct sum of all indecomposable
projective-injective modules in $\Oo_0^s$ and consider $\PT$ as an
object in $\Oo_0$. We claim that $\Phi: \varphi\mapsto \varphi_{Q}$
defines an isomorphism $\End(\ZZ)\cong \cZ(\End_{\g}(\PT))$, where
the latter denotes the center of $\END_\mg(\PT)$. Note that
$\cZ(\End_\g(\PT))\cong \END(\ID_{\cO_0^s})$ (\cite[Theorem
10.1]{Stquiv}).

{\it $\Phi$ is injective:} Let $\varphi\in\END(\ZZ)$,
$\varphi_{\PT}=0$ and  let  $P$ be a projective object in $\Oo_0$.
We fix an inclusion $i:\ZZ P\inj J_1$, where $J_1=\oplus_{i\in
I_1}\PT$ for some finite set $I_1$ (see the main result of \cite{I2}). Since $\ZZ$ is
the identity on $\cO_0^s$ we have $\phi_P=\phi_{ZP}$ and
$0=\phi_{J_1}\circ \ZZ(i)=\ZZ(i)\circ\phi_{\ZZ P}$. The injectivity
of $\ZZ(i)$ implies $\varphi_P=0$. Let $M\in\cO_0$ be arbitrary with
projective cover $f:P\surj M$. Then $\varphi_M\circ \ZZ(f)=
\ZZ(f)\circ\phi_{J_i}$, i.e. $\varphi_M=0$, since $\ZZ$ is right
exact.

{\it $\Phi$ is surjective:} Let $g\in\cZ(\End_{\g}(\PT))$. For $P\in
\Oo_0$ projective we fix a coresolution
\begin{eqnarray*}
  \ZZ P\stackrel{i}\inj J_1\stackrel{h}\longrightarrow J_2,
\end{eqnarray*}
where $J_i\cong\oplus_{i\in I_i}\PT$ for some finite sets $I_i$
($i=1,2$). For the existence of such a tilting resolution one can use the main result of
\cite{I2} and arguments, analogous to that of \cite[3.1]{KSX} (see
\cite[Theorem~10.1]{Stquiv}).
By definition, $g$ induces a natural map
$g_{\ZZ P}\in\END_\mg(\ZZ P)$ making the following
diagram commutative:
\begin{displaymath}
\xymatrix{
\ZZ P \ar@{^{(}->}[rr]^{\ZZ(f)}\ar@{.>}[d]_{g_P}&&
\ZZ J_1\ar[rr]^{\ZZ(h)}\ar[d]_{g_{J_1}}
&& \ZZ J_2\ar[d]_{g_{J_2}}\\
\ZZ P \ar@{^{(}->}[rr]^{\ZZ(f)}&& \ZZ J_1\ar[rr]^{\ZZ(h)} && \ZZ J_2\\
}.
\end{displaymath}
Setting $g_P=g_{\ZZ P}$ defines a natural transformation
$\tilde{g}:\ZZ\rightarrow\ZZ$, when restricted to the additive
category of projective objects in $\cO_0$ such that
$\tilde{g}_{\PT}=g$. The right exactness of $\ZZ$ ensures that
$\tilde{g}$ extends uniquely to some $\tilde{g}\in\END(\ZZ)$.
Hence $\Phi$ is surjective. In particular,
$\End(\ZZ)=\cZ(\End_{\g}(\PT))=\cZ(\Oo_0^s)=\END(\ID_{\Oo_0^s})$.
\end{proof}

The remaining part from Theorem~\ref{tmain4} follows if we
prove the following statements:

\begin{proposition}\label{pr6.20}
Let $F:\mathcal{A}\to \mathcal{B}$ be a dense functor between two
categories $\mathcal{A}$ and $\mathcal{B}$. Then the restriction
gives rise to an injective linear map
$\END(\ID_\mathcal{B})\hookrightarrow\End(F)$. In particular,
$\ZZ\FF:\cO_0\rightarrow\cO_0^s$ provides an inclusion
$\END(\ID_{\cO_0^s})\hookrightarrow\ho_{\GG,\TT}$.
\end{proposition}

\begin{proof}
The first statement of the proposition is obvious. Since
$\ZZ\FF M =M$ for any $M\in\cO_0$ we may consider $\FF=\ZZ\FF$ as
a functor from $\cO_0$ to $\cO_0^s$. We claim that $\FF$ is dense,
i.e. for any $N\in \Oo_0^s$ there exists an $K\in \Oo_0$ such that
$\ZZ\FF(K)\cong N$. Indeed, let $P\tto N$ be a projective cover of $N$
in $\Oo_0$ with kernel $K$. Applying $\GG$ to $K\hookrightarrow P \tto
N$ we obtain the exact sequence $\GG K\hookrightarrow \GG P \to \GG N$
and $\GG N=0$. In particular, $\GG K\cong \GG P$. Since the socle of
$P$, and therefore also of $K$, is annihilated by $\ZZ$, the map
$\rg_K$ is injective (see \cite[Lemma~2.4]{J}). Hence we have $\FF K\cong
(\GG K)/ K\cong (\GG P)/K\cong P/K\cong N$.

By Theorem~\ref{tmain1} we have morphisms
$\GG\overset{\rp}{\longrightarrow}\FF
\overset{\alpha^{-1}}{\longrightarrow}\FF'
\overset{\ri}{\hookrightarrow}\TT$,
where $\alpha^{-1}$ is an isomorphism. We consider the linear map
$\xi:\End(\FF)\to \ho_{\GG,\TT}$ defined as
$\xi(\varphi)=\ri\circ \alpha^{-1}\circ
\varphi\circ\rp$. Since $\rp$ is surjective, $\ri$ is injective, and
$\alpha^{-1}$ is an isomorphism, $\xi$ defines an inclusion
$\End(\FF)\hookrightarrow \ho_{\GG,\TT}$. To complete the proof
it is now enough to show that $\End(\FF)$ contains $\END(\ID_{\Oo_0^s})$. This
follows directly from the first part of the proposition, since
$\END(\ZZ)\cong \END(\ID_{\Oo_0^s})$ (by Proposition~\ref{l6.25}).
\end{proof}

\begin{remark}
{\rm The case $\mg=\mathfrak{sl}_2$ shows already that some spaces
$\ho_{X,Y}$, $X,Y\in\mathcal{S}$ can be smaller than $\Coinv$. Indeed,
in this case we have $\ho_{\GG,\ID}\cong\mC$ and
$\ho_{\GG\TT,\TT\GG}\cong\mC$. Although
the remaining `unknown' spaces from Theorem~\ref{tmain4} are
isomorphic to $\Coinv$ in this particular example, the isomorphism is
accidental and is not given by a natural action of $\Coinv$ on
$P(w_0)$ (in contrast to the cases, which are known to be isomorphic
to $\Coinv$ from Theorem~\ref{tmain4}).
\erem}
\end{remark}

\section{Proof of Theorem~\ref{tmain4.101}}\label{s6.101}

Let $\mathcal{I}(\mathcal{\hat{S}})$ denote the set of all idempotents
in $\mathcal{\hat{S}}$. For $\XX,\YY\in\mathcal{I}(\mathcal{\hat{S}})$
we set $\ho_{\XX,\YY}=\Hom(\XX,\YY)$.

\begin{proposition}\label{pr6.105}
$\End(\XX)\cong \Coinv$ as algebras for any
$\XX\in\mathcal{\hat{S}}$.
\end{proposition}

\begin{proof}
An injective algebra morphism from $\Coinv$ to $\End(\XX)$ for every
$\XX\in\mathcal{\hat{S}}$ is constructed using the same arguments as
in Proposition~\ref{pr6.1}. The arguments, analogous to that of
Proposition~\ref{pr6.1}, also give an isomorphism
$\End(\CC)\cong\Coinv$.

Let us show that $\End(\CC^2)\cong\Coinv$. We claim that the
evaluation $\phi\mapsto\phi_{P(w_0)}$ defines an inclusion
$\END(\CC^2)\inj\END_\mg(\CC^2P(w_0\cdot0)$. Assume $\phi_{P(w_0)}=0$
and let $P\in\cO_0$ be projective with injective hull $i:P\inj I$.
We get an exact sequence
$0\rightarrow\KER\CC^2(i)\rightarrow\CC^2P\rightarrow\CC^2 I$. By
assumption we have $0=\phi_I\circ\CC^2(i)=\CC^2(i)\circ\phi_P$. In
particular, the image of $\phi_P$ is contained in the kernel of
$\CC^2(i)$. On the other hand $\HOM_\mg(\CC^2
P,\KER\CC^2(i))\inj\HOM_\mg(\theta\CC P,\KER\CC^2(i))\cong\HOM(\CC
P,\theta\KER\CC^2(i))=0$, since $\theta\KER\CC^2(i)=0$. Therefore,
$\phi_P=0$ and hence $\phi=0$, since $\CC^2$ is right exact.

If $i>2$ then we have
\begin{displaymath}
\End(\CC^i)\cong\Hom(\ID,\KK^i\CC^i)\cong\Hom(\ID,\KK^2\CC^2)
\cong\End(\CC^2)\cong \Coinv.
\end{displaymath}
$\End(\KK\CC^i)\cong
\Hom(\CC\KK\CC^i,\CC^i)\cong \End(\CC^i)\cong\Coinv$, $i>0$; and
$\End(\KK^2\CC^i)\cong
\Hom(\CC^2\KK^2\CC^i,\CC^i)\cong \End(\CC^i)\cong\Coinv$, $i>1$.

Finally, there are isomorphisms
\begin{multline*}
\End(\CC\KK^2\CC)\cong\Hom(\KK^2\CC,\KK\CC\KK^2\CC)\cong
\End(\KK^2\CC)\cong\\ \cong  \Hom(\CC^2\KK^2\CC,\CC)\cong
\Hom(\CC^2\KK,\CC)\cong\Hom(\CC\KK,\KK\CC)
\end{multline*}
and it is left to show
that $\Hom(\CC\KK,\KK\CC)$ embeds into $\Coinv$ as a vector
space. For this we show that the map $\Phi:\Hom(\CC\KK,\KK\CC)\to
\End_{\mathfrak{g}}(P(w_0))\cong \Coinv$,
$\varphi\mapsto\varphi_{P(w_0)}$ is injective. Assume that
$\varphi_{P(w_0)}=0$. Since both $\CC\KK$ and $\KK\CC$
preserve injections (see Proposition~\ref{l2.102}), from the injection
$i:P\inj I$ above we obtain that
$\varphi$ must be zero on all projective modules. Taking a projective
cover of any $M\in\Oo_0$ and using the fact that both $\CC\KK$ and
$\KK\CC$ preserve surjections (see Proposition~\ref{l2.102}), we obtain
that $\varphi$ is zero.
The rest follows by duality.
\end{proof}

Note that $\KK\CC$ preserves projective modules, since the adjunction from
Proposition~\ref{l2.102} is an isomorphism on projective objects.

{\it Equality of the spaces labeled by $2$}: The inclusion
$\CC\KK\hookrightarrow \ID$ from Proposition~\ref{l2.102} induces an inclusion
$\ho_{\KK^2\CC^2,\CC\KK}\hookrightarrow \ho_{\KK^2\CC^2,\ID}$.
By duality we have $\ho_{\KK^2\CC^2,\CC\KK}\cong
\ho_{\KK\CC,\CC^2\KK^2}$ and $\ho_{\KK^2\CC^2,\ID}\cong
\ho_{\ID,\CC^2\KK^2}$. Applying  Proposition~\ref{p6.4} to
$F=\ID$, $H=\KK\CC$, and $J=\CC^2\KK^2$ we obtain
$\ho_{\ID,\CC^2\KK^2}\hookrightarrow \ho_{\KK\CC,\CC^2\KK^2}$ and thus
all these four spaces are isomorphic.

{\it Equality of the spaces labeled by $3$}: The inclusion
$\CC\KK\hookrightarrow \ID$ induces the following inclusion:
$\ho_{\KK\CC^2\KK,\CC\KK}\hookrightarrow \ho_{\KK\CC^2\KK,\ID}$.
By duality we have $\ho_{\KK\CC^2\KK,\CC\KK}\cong
\ho_{\KK\CC,\CC\KK^2\CC}$ and $\ho_{\KK\CC^2\KK,\ID}\cong
\ho_{\ID,\CC\KK^2\CC}$. Applying  Proposition~\ref{p6.4} to
$F=\ID$, $H=\KK\CC$, and $J=\CC\KK^2\CC$ we obtain
$\ho_{\ID,\CC\KK^2\CC}\hookrightarrow \ho_{\KK\CC,\CC\KK^2\CC}$ and
thus all these four spaces are isomorphic.

{\it Equality of the spaces labeled by $4$}:
Evaluating $\CC\KK\hookrightarrow\ID$ at $\KK\CC$ gives an inclusion
$\CC\KK^2\CC\cong \KK\CC^2\KK\hookrightarrow\KK\CC$.
Applying $\Hom(\KK^2\CC^2,{}_-)$
produces $\ho_{\KK^2\CC^2,\KK\CC^2\KK}\hookrightarrow
\ho_{\KK^2\CC^2,\KK\CC}$. By duality we have
$\ho_{\KK^2\CC^2,\KK\CC^2\KK}\cong
\ho_{\KK\CC^2\KK,\CC^2\KK^2}$ and $\ho_{\KK^2\CC^2,\KK\CC}\cong
\ho_{\CC\KK,\CC^2\KK^2}$. Applying Proposition~\ref{p6.4} to
$F=\CC\KK$, $H=\KK\CC$, and  $J=\CC^2\KK^2$ we obtain
$\ho_{\CC\KK,\CC^2\KK^2}\hookrightarrow
\ho_{\CC\KK^2\CC,\CC^2\KK^2}$ and thus all these four spaces are
isomorphic.

Applying the duality implies that all other spaces labeled by
the same number coincide.

{\it All spaces labeled by  $\Coinv$ are correct}:
For the diagonal entries this follows from Proposition~\ref{pr6.105}
above. For any $\XX\in \mathcal{I}(\mathcal{\hat{S}})$ we have
$\ho_{\CC^2\KK^2,\XX}\cong \ho_{\KK^2,\KK^2\XX}\cong
\ho_{\KK^2,\KK^2}\cong \Coinv$ and $\ho_{\XX,\KK^2\CC^2}\cong\Coinv$
by duality. That $\ho_{\CC\KK,\KK\CC}\cong\Coinv$ was shown in
the proof of Proposition~\ref{pr6.105}. Using adjunction and duality we have
$\ho_{\CC\KK,\KK\CC^2\KK}\cong\ho_{\CC^2\KK,\CC^2\KK}\cong\Coinv$ and $\ho_{\ID,\KK\CC}\cong\ho_{\CC,\CC}\cong\Coinv\cong\ho_{\KK,\KK}\cong\ho_{\CC\KK,\ID}$.

It is left to establish the claimed inclusions. Applying
$\Hom(\KK\CC,{}_-)$ to the inclusion $\CC\KK\hookrightarrow\ID$ we
get $\ho_{\KK\CC,\CC\KK}\hookrightarrow\ho_{\KK\CC,\ID}$. Applying
Proposition~\ref{p6.4} to $F=\ID$, $H=\KK\CC$, and $J=\CC\KK$ we
obtain $\ho_{\ID,\CC\KK}\hookrightarrow\ho_{\KK\CC,\CC\KK^2\CC}$.
Applying Proposition~\ref{p6.4} to $F=\KK\CC^2\KK$, $H=\KK\CC$, and
$J=\CC\KK$ we obtain $\ho_{\KK\CC^2\KK,\CC\KK}\hookrightarrow
\ho_{\CC\KK^2\CC,\CC\KK^2\CC}\cong\Coinv$.
Applying $\Hom(\KK^2\CC^2,{}_-)$ to the inclusion
$\KK\CC\hookrightarrow\KK^2\CC^2$ obtained above we get
$\ho_{\KK^2\CC^2,\KK\CC}\hookrightarrow
\ho_{\KK^2\CC^2,\KK^2\CC^2}\cong\Coinv$.

\begin{remark}{\rm
Behind our argumentation is the following general fact:
Let $F$ and $G$ be two endofunctors on $\Oo_0$. Assume that
$F$ preserves surjections and $G$ preserves injections.
Then the map $\Hom(F,G)\to \End_{\mathfrak{g}}(P(w_0))$,
$\varphi\mapsto \varphi_{P(w_0)}$, is injective. Indeed,
let $\varphi_{P(w_0)}=0$. Since the injective envelope of any
projective $P\in \Oo_0$ belongs to $\mathrm{Add}(P(w_0))$, we can use
that $G$ preserves injections to obtain $\varphi_{P}=0$. Taking
now the projective cover of any $M\in\Oo_0$ and using that $F$
preserves surjections we obtain $\varphi_{M}=0$.

One can show that $\KK^2\CC^2$ preserves injections and $\CC^2\KK^2$
preserves surjections, which implies that
$\ho_{\XX,\YY}\hookrightarrow\Coinv$ for all
$\XX\in\{\ID,\CC\KK,\KK\CC^2\KK,\CC^2\KK^2,\CC^i,\KK\CC^i:i>0\}$
and for all $\YY\in\{\ID,\KK\CC,\KK\CC^2\KK,\KK^2\CC^2,
\KK^i,\CC\KK^i:i>0\}$.
\erem}
\end{remark}

\section{Proof of Theorem~\ref{tmain5}}\label{s7}

We have  $\Ext_{\Oo_0}^{i}(\mathcal{P}^w,\mathcal{P}^w)=
\Hom_{D^b(\cO_0)}(\cL\TT_w\mathcal{P}, \cL\TT_w\mathcal{P}[i])=
\EXT^{i}_{\cO_0}(\mathcal{P},\mathcal{P})=0$, $i>0$, (see
\cite[Corollary~4.2]{AS}).

\begin{claim}\label{cl7.2}
$\mathcal{P}$ admits a finite coresolution by
modules from $\mathrm{Add}(\mathcal{P}^{w})$.
\end{claim}

\begin{proof}
Let $w\in W$. If $l(w)=0$, the statement is obvious. Assume, it is
true for all $\tilde{w}$ where $l(\tilde{w})\leq l(w)$ and let $s$
be a simple reflection such that $sw>w$. We have to show that
$\mathcal{P}$ has a finite
$\mathrm{Add}(\mathcal{P}^{sw})$-coresolution.
Since $\Ext_{\Oo_0}^{>0}(\mathcal{P}^x,\mathcal{P}^x)=0$, for all
$x\in W$, the arguments from \cite[Chapter~III]{Ha} or
\cite[Lemma~4]{MO} reduce the problem to showing that there exists a
$\tilde{w}$, $l(\tilde{w})\leq l(w)$, such that $P^{\tilde{w}}$ admits
a coresolution by modules from  $\mathrm{Add}(\mathcal{P}^{sw})$.
Since all $\TT_x$ commute with translation functors, it is enough to
prove the statement for $\TT_{\hat{w}}\Delta(e)\cong\Delta(\hat{w})$.
We choose $\tilde{w}$ such that $sw=\tilde{w}t$ for some simple
reflection $t$ with $l(\tilde{w}t)>l(\tilde{w})$ and consider the
short exact sequence $\Delta(e)\hookrightarrow P(t)\tto\Delta(t)$.
Applying $\TT_{\tilde{w}}$ we obtain the short exact sequence
$\Delta(\tilde{w})\hookrightarrow \TT_{\tilde{w}}
P(t)\tto\Delta(sw)$. Since $P(t)\cong \TT_t P(t)$, it follows that
$\TT_{\tilde{w}} P(t)\cong \TT_{sw}P(t)$. Thus,
$\TT_{\tilde{w}} P(t), \Delta(sw)\in \mathrm{Add}(\mathcal{P}^{sw})$,
and hence $\Delta(\tilde{w})$ has a coresolution by modules from
$\mathrm{Add}(\mathcal{P}^{sw})$.
\end{proof}

We proved that $\mathcal{P}^w$ is a generalized tilting module for
any $w\in W$. Since $\cO_0$ has finite projective dimension, it is
a generalized cotilting module as well (\cite[Corollary~2.4]{Reiten}).

The remaining assertions from the first part of the theorem follow by
duality. Since $\TT_{w_0}\Delta(e)\cong\Delta(w_0)$ is a tilting
module and $\TT_{w_0}$ commutes with translations, it follows that
$\mathcal{P}^{w_0}\cong \mathcal{T}\cong\mathcal{I}^{w_0}$.
Let $w\in W$ and $sw>w$ (i.e. $sww_0<ww_0$). The
adjunction morphism $\TT_s\GG_s\hookrightarrow \ID$ gives
$\TT_{sw}\TT_{w_0}\mathcal{P}\cong\TT_s\TT_w\TT_{w_0}\mathcal{P}
\cong\TT_s\GG_{ww_0}\mathcal{I}\cong
\TT_s\GG_s\GG_{sww_0}\mathcal{I}\hookrightarrow
\GG_{sww_0}\mathcal{I}$. Comparing the characters and using duality
shows the second part of the theorem.

It remains to prove the formulas for the homological dimensions.
Twisting functors commute with translation functors, hence we get
\begin{displaymath}
\mathrm{projdim}(\mathcal{P}^{w})=\mathrm{projdim}(\TT_w\Delta(e))=
\mathrm{projdim}(\Delta(w))
\end{displaymath}
and $\mathrm{injdim}(\mathcal{P}^{w})=
\mathrm{injdim}(\Delta(w))$. For Verma modules the values are easy to
compute and are given by the formulas from the theorem. The remaining
statements follow by duality.

\section{Proof of Theorem~\ref{tmain6}}\label{s8}

We start with the following

\begin{proposition}\label{derived}
Let $w\in W$ and $M\in\cO_0$ be a module having a Verma flag. Then
$\cL_1\CC_s(\CC_{w^{-1}} M)=0$ for any simple reflection $s$ such that
$ws>w$. In particular, $\CC_{w^{-1}} P$ is acyclic for $\CC_s$ for
any projective object $P$ and hence $\cL\CC_s\cL\CC_{w^{-1}}\cong
\cL(\CC_s\CC_w)$.
\end{proposition}

\begin{proof}
By \cite[Section~5]{MS}, $\CC_{w^{-1}} M$ has a $w^{-1}$-shuffled
Verma flag. Hence, using Theorem~\ref{tmain2}, it is enough to show
that the socle of every $w^{-1}$-shuffled Verma module
$\CC_{w^{-1}} \Delta(x)$ contains only $L(y)$ such that $ys<y$.
But $\CC_{w^{-1}} \Delta(x)$ is at the same time
a $w^{-1}w_0$-coshuffled dual Verma module and $sw^{-1}w_0<w^{-1}w_0$
as $ws>w$. This implies that $\CC_{w^{-1}} \Delta(x)\cong \KK_{s} N$
for some $N\in \Oo_0$ and thus $\CC_{w^{-1}} \Delta(x)$ has desired
socle by definition of $\KK_s$.
\end{proof}

\begin{claim}
${}^w\mathcal{P}$ is a generalized (co-)tilting module.
\end{claim}

\begin{proof}
The case $w=e$ is clear. Assume the statement to be true for $w\in W$
and let $s$ be a simple reflection such that $sw>w$. By definition
\begin{eqnarray*}
0\rightarrow P(x)\stackrel{\ADJ_s(P(x))}\longrightarrow
\theta_sP(x)\longrightarrow \CC_sP(x)\rightarrow 0
\end{eqnarray*}
is exact for any $x\in W$. Applying $\CC_w$ and using the previous
proposition we get an exact sequence
\begin{eqnarray*}
0\rightarrow\CC_w P(x)\stackrel{\ADJ_s(P(x))}\longrightarrow
\CC_w\theta_sP(x)\longrightarrow \CC_w\CC_s P(x)\rightarrow 0.
\end{eqnarray*}
Since $\CC_w\CC_s\cong\CC_{sw}$ (see \cite[Lemma~5.10]{MS}) and
$\CC_w\theta_sP(x)\cong\CC_w\CC_s\theta_sP(x)\cong\CC_{sw}\theta_sP(x)$,
$\CC_wP(x)$ has a two-step coresolution with modules from
$\mathrm{Add}(\CC_{ws}\mathcal{P})$. Since $\cL\CC_w$ induces an
equivalence on the bounded derived category of $\cO_0$ (by
Proposition~\ref{derived} and \cite[Theorem~5.7]{MS}) we have
$\EXT^{>0}(\CC_w\mathcal{P},\CC_w\mathcal{P})\cong
\EXT^{>0}(\mathcal{P},\mathcal{P})$.
The arguments from Claim~\ref{cl7.2} show that $_{}^w\cP$ is a
generalized tilting module, hence also a generalized cotilting module
by \cite[Corollary~2.4]{Reiten}.
\end{proof}

Now let us  prove Theorem~\ref{tmain6}\eqref{tmain6.3}. Using
Proposition~\ref{derived} and \cite[Section~5]{MS} the statement
reduces to verifying that $_{}^{w_0}\mathcal{P}\cong
\mathcal{T}$. Since $\CC_{w_0}$ maps Verma modules to dual Verma
modules, Proposition~\ref{derived} implies that
$\CC_{w_0}\mathcal{P}$ has a dual Verma flag and satisfies
$\Ext_{\Oo_0}^{i}\left(\CC_{w_0}\mathcal{P},\rd\,\Delta(x)\right)=0$
for all $x\in W$. From \cite[Corollary~4]{R} it follows that
$\CC_{w_0}\mathcal{P}$ has a Verma flag as well and thus
$\CC_{w_0}\mathcal{P}\cong \mathcal{T}$.

Let $L=L(y)\in\cO_0$ be a simple object and $M\in\Oo_0$ be a module
with Verma flag. Then Proposition~\ref{derived} gives
\begin{eqnarray*}
\EXT^i_{\cO}(\CC_s\CC_w M,L)&\cong&
\HOM_{D^b(\cO_0)}(\cL(\CC_s\CC_w)M,L[i])\\
&\cong&
\HOM_{D^b(\cO_0)}(\CC_w M,\cR\KK_s
L[i]).
\end{eqnarray*}
The latter is $\EXT^{i+1}_{\cO}(\CC_w M,L)$ if $y<ys$
and it is $\EXT^{i}_{\cO}(\CC_w M,K_s L)$ otherwise (see
\cite[Proposition~5.3]{MS}). In particular, $M=\mathcal{P}$ gives
$\mathrm{projdim}({}^{ws}\mathcal{P})\leq
\mathrm{projdim}({}^{w}\mathcal{P})+1$, and $M=\mathcal{T}$ gives
$\mathrm{projdim}(C_{ws}\mathcal{T})\leq
\mathrm{projdim}(C_{w}\mathcal{T})+1$. However, we know that
$\mathrm{projdim}(\mathcal{T})=\mathrm{injdim}(\mathcal{T})=l(w_0)$
(see e.g. \cite[\S~7]{BGG}, \cite[Theorem~1]{MO}) and $\mathrm{projdim}(\mathcal{I})=l(w_0)$
and all the formulae for homological dimensions follow.

\begin{remark}
{\rm It is well-known (see e.g. \cite[Theorem~2.1 and Section~3]{AL}) that
the set of twisted  Verma modules is equal to the set of shuffled Verma
modules. This is not the case for projective objects. In fact, if
$\mg=\SL_3$ and $s$, $t$ are the two simple reflections, then
direct calculations show that $\CC_sP(t)$ is neither a twisted
projective nor a completed injective object.
\erem}
\end{remark}

\section{Proof of Theorem~\ref{singularbraid}}
\label{BraidMonoid}
The {\it singular braid monoid} is generated by
    $\{\sigma_i,\sigma_i^{-1},\tau_i\}$ ($1\leq i\leq n-1$) subject to the
    relations
    \begin{eqnarray}
\label{B1}
   \sigma_i\sigma_i^{-1}=\sigma_i^{-1}\sigma_i&=&1,\mbox{    for all } i,\\
\label{B2}
        \sigma_i\sigma_{i+1}\sigma_i&=&\sigma_{i+1}\sigma_i\sigma_{i+1},\mbox{   for all } i,\\
\label{B3}
      \sigma_i\sigma_j&=&\sigma_j\sigma_i\mbox{   if } |i-j|>1,\\
\label{B4}
\tau_{i}\sigma_j\sigma_{i}&=&\sigma_j\sigma_{i}\tau_j\mbox{    if   } |i-j|=1,\\\label{B6}
\sigma_i\tau_j&=&\tau_j\sigma_i\mbox{    if }|i-j|\not=1,\\
\label{B7}
\tau_i\tau_j&=&\tau_j\tau_i\mbox{    if }|i-j|>1.
\end{eqnarray}

For a different presentation see for example \cite{DasbachGemein}.
We have to prove that the functors from Theorem~\ref{singularbraid} satisfy
the relations. For the first three relations see e.g. \cite[Theorem~5.7 and Lemma~5.10]{MS}.
We claim that the remaining relations are true on the level of
endofunctors on $\cO_0$. Then they are also true for the derived functors (note
that $\cL(\CC_s\CC_t)\cong\cL\CC_s\cL\CC_t$ if $s\not=t$ by e.g. \cite[Proposition~3.1]{I}).
Relation~\eqref{B7} follows directly from the classification
theorem (\cite[3.3]{BG}) of projective functors. If $i\not=j$, then the relation
\eqref{B6} follows immediately from the definition of $\CC_s$. In the case $i=j$
the relation \eqref{B6} will be proved in Lemma~\ref{thetaC} below. The relations
\eqref{B4} will be proved in  Proposition~\ref{difficult}.

\begin{lemma}
  \label{thetaC}
With the notation from Theorem~\ref{singularbraid} we have: There are isomorphisms of functors
$\CC_{s_i}\theta_{s_i}\cong\theta_{s_i}\CC_{s_i}$ for $1\leq i\leq n-1$.
\end{lemma}

\begin{proof}
  We set $s=s_i$ for some $1\leq i\leq n-1$.
  All occurring functors are right exact and exact on modules having a Verma
  flag (see \cite[Proposition 5.3]{MS}). Note that they preserve the
  full subcategory $\cT$ of projective-injective
  modules in $\cO_0$. We claim that it is enough to establish the isomorphism when
  restricted to this category. Indeed, any projective object has a coresolution
  by objects in $\cT$, then standard arguments using the Five Lemma will extend
  the constructed
  isomorphism to an isomorphism of the corresponding endofunctors on the
  category of projective modules  (since $\theta_s \CC_s M\cong \theta_s M$
  for any object $M$, all functors in question preserve this category).
  Again by standard arguments, using projective resolutions, the statement
  would follow, since the functors are right exact.

Hence, let us consider the category $\cT$. The functor $\mV$
  from Section~\ref{s2} defines an equivalence $\tilde{\mV}$ of categories between $\cT$ and
  the category of finite dimensional free $\cC$-modules. We have
  $\tilde\mV\theta_s\tilde\mV^{-1}$ is given by tensoring with the bimodule
  $\cC\otimes_{\cC^s} \cC$ (see \cite[Theorem~10]{S}). Recall that $\cC$ is a free
  $\cC^s$-module of rank $2$ with basis $1$ and
  $X$, the coroot corresponding to $s$. From the definitions it follows then
  that $\CC_s\theta_s$ is given by tensoring with the cokernel $D_s$ of the
  map
  \begin{eqnarray*}
    \phi:\quad\cC\otimes_{\cC^s}\cC&\stackrel{\op{adj_s}\otimes\op{id}}{\rightarrow}&
    \cC\otimes_{\cC^s} \cC\otimes_{\cC^s}\cC=:E\\
    c\otimes d&\mapsto& X\otimes c\otimes d+1\otimes cX\otimes d,
  \end{eqnarray*}
where $\op{\op{adj}_s}$ denotes the corresponding adjunction map.

We define a homomorphism of vector spaces from $E$ to $\cC\otimes_{\cC^{s}}\cC$ by
\begin{displaymath}
\begin{array}{ll}
1\otimes 1\otimes d \mapsto 1\otimes d, & X\otimes 1\otimes d \mapsto
 X\otimes d, \\ 1\otimes X\otimes d \mapsto -X\otimes d, & X\otimes X\otimes d
 \mapsto -X^2\otimes d
\end{array}
\end{displaymath}
for any $d\in \cC$. This is obviously well-defined and
 defines in fact a unique $\mC$-linear map. Evidently, it factors through $D_s$, is surjective, and is a homomorphism of
 $\cC$-bimodules. Since $\CC_s\theta_s P(w_0)\cong\theta_s P(w_0)$ it has to be an
 isomorphism. Hence $\CC_s\theta_s\cong \theta_s$ on $\cT$ for any simple
 reflection $s$. By the remarks above we get  $\CC_s\theta_s\cong \theta_s$ as
 endofunctors on $\cO_0$.
 Similarly one proves that $\theta_s\CC_s\cong \theta_s$ by looking at the
 cokernel of the map
   \begin{eqnarray*}
    \cC\otimes_{\cC^s}\cC&\stackrel{\op{id}\otimes\op{adj_s}}{\rightarrow}& \cC\otimes_{\cC^s} \cC\otimes_{\cC^s}\cC\\
    c\otimes d&\mapsto& c\otimes X\otimes d+c\otimes 1\otimes Xd.
  \end{eqnarray*}
The statement of the lemma follows.
\end{proof}

\begin{remark}
 {\rm Using the graded version from \cite[Section~7]{MS} the isomorphism from the
  previous lemma is given as
  follows: We choose an isomorphism of functors $\phi:\theta_s\langle
  1\rangle\oplus\theta_s\langle
  -1\rangle\cong\theta_s^2$. The isomorphism $\theta_s\langle -1\rangle\cong\theta_s\CC_s$ is
  then given as $\theta_sp\circ\phi\circ i_2$, where $p$ is the canonical
  projection and $i_2$ denotes the inclusion into the second summand.
\erem}
\end{remark}

  \begin{proposition}
\label{difficult}
    With the notation and assumptions from Theorem~\ref{singularbraid} we have isomorphisms
    of functors
    \begin{eqnarray*}
    \theta_{i+1}\CC_i\CC_{i+1}&\cong&\CC_i\CC_{i+1}\theta_i\mbox{    for
      all }i,\\
    \CC_{i+1}\CC_{i}\theta_{i+1}&\cong&\theta_{i}\CC_{i+1}\CC_{i} \mbox{
      for all }i.
    \end{eqnarray*}
 \end{proposition}

 \begin{proof}
   We give an argument for the first
   isomorphism, and omit the analogous calculations for the second one. Set $s=s_{i+1}$,
   $t=s_i$.
Note first that it is sufficient to establish the isomorphism on projective
modules, since the functors are right exact. Since any projective module has a
   copresentation by projective-injective modules and since the functors in
   question are exact on modules with Verma flag (\cite[Proposition~3.1]{I} and \cite[Proposition~5.3]{MS}), it is enough to
   check it on the subcategory given by these objects. We first compare them
   on the Grothendieck group level. Since the functors in question are exact on modules with Verma flag, we may even restrict ourselves to the case of a single Verma
module. We have the following formulas in the Grothendieck group of $\cO_0$:
$[\theta_s\CC_t\CC_s\Delta(x)]=[\TT_{x}\theta_s\CC_t\CC_s\Delta(e)]=[\TT_x\theta_s\Delta(st)]=[\TT_{x}(\Delta(st)\oplus\Delta(sts))]=[\Delta(xst)\oplus\Delta(xsts)]$.
    Here $\TT_{x}=\TT_{s_1}\TT_{s_2}\cdots\TT_{s_r}$, where $x=s_1s_2\cdots
    s_r$ is a reduced expression (see e.g. \cite[Section~6]{KM}). On the other hand we
    have $[\CC_t\CC_s\theta_t\Delta(x)]=[\TT_{x}\CC_t\CC_s\theta_t\Delta(e)]=[\TT_x\CC_t\CC_s(\Delta(e)\oplus\Delta(t)]=[\TT_{x}(\Delta(st)\oplus\Delta(tst))]=[\Delta(xst)\oplus\Delta(xsts)]$.

Before we proceed, we want to give the principal
   idea of the proof. The classification theorem of projective functors
   (\cite[3.3]{BG})
   provides (in the case of $\mg=\SL_n$) a decomposition $\theta_s\theta_t\theta_s\cong F\oplus\theta_s$
   for any noncommuting simple reflections $s$ and $t$. Here, $F$ is the indecomposable functor given by
   $F(\Delta(e))=P(sts)$. By the definition of the functors we get surjections
   $\alpha$ and $\beta$ defined by the following commuting diagrams:
   \begin{eqnarray*}
\xymatrix{
\theta_s\theta_t\theta_s\ar@{->>}[rd]|{\alpha}\ar@{->>}[r]^{\alpha_1:=\theta_s\theta_t
     p}\ar@{->>}[d]_{\alpha_2:=\theta_s(p_{\CC_s})}
&\theta_s\theta_t\CC_s\ar@{->>}[d]^{\theta_s
     (p_{\theta_s})}\\
\theta_s\CC_t\theta_s\ar@{->>}[r]_{\theta_s\CC_t(p)}&\theta_s\CC_t\CC_s,}
&&
\xymatrix{
\theta_t\theta_s\theta_t\ar@{->>}[rd]|{\beta}\ar@{->>}[r]^{\beta_1:=\theta_tp_{\theta_t}}
\ar@{->>}[d]_{\beta_2:=p_{\theta_s\theta_t}}
&\theta_t\CC_s\theta_t\ar@{->>}[d]^{p_{\CC_s\theta_t}}\\
\CC_t\theta_s\theta_t\ar@{->>}_{\CC_t(p_{\theta_t})}[r]&\CC_t\CC_s\theta_s,
}
\end{eqnarray*}
where $p$ always denotes the corresponding natural projection.
We claim that already $\alpha\circ i$ is surjective and factors through
$\beta\circ j$ for some fixed inclusions
$i:F\rightarrow\theta_s\theta_t\theta_s$ and
$j:F\rightarrow\theta_t\theta_s\theta_t$; i.e. there exists a natural transformation
$h:\theta_s\CC_t\CC_s\rightarrow \CC_t\CC_s\theta_t$ which is a
surjection. The statement would then follow from our comparison on the
Grothendieck group level.

As in the proof of the previous lemma we will work with $\cC$-bimodules. The map
$\alpha$ gives then rise to an endomorphism of $\cC$-bimodules
\begin{eqnarray*}
 \tilde{\alpha}: &&\cC\otimes_{\cC^s}\cC\otimes_{\cC^t} \cC\otimes_{\cC^s} \cC\rightarrow D_\alpha,\\
 \tilde{\beta}: &&\cC\otimes_{\cC^t}\cC\otimes_{\cC^s} \cC\otimes_{\cC^t} \cC\rightarrow D_\beta,
\end{eqnarray*}
where $D_\alpha\otimes_\cC\bullet\cong\mV\theta_s\CC_t\CC_s\mV^{-1}$ and
$D_\beta\otimes_\cC\bullet\cong\mV\CC_t\CC_s\theta_t\mV^{-1}$ on the category of
free $\cC$-modules of finite rank (see the proof of
the previous lemma). Set $\tilde\alpha_i=\mV\alpha_i\mV^{-1}$ and
$\tilde\beta_i=\mV\beta_i\mV^{-1}$ for $i=1$, $
2$.
Let $X$ and
$Y$ be
the coroots corresponding to $s$ and $t$ respectively. Note that the $\{b\otimes y\otimes x\otimes c\}$  for $c$ running through
a $\mC$-basis of $\cC$, and $b$, $x\in\{1,X\}$, $y\in\{1,Y\}$ are a $\mC$ basis
of $\cC\otimes_{\cC^s}\cC\otimes_{\cC^t} \cC\otimes_{\cC^s} \cC$ (this
follows from the fact that $\cC$ is a free $\cC^s$ module with basis $1$ and $X$).\\
{\it We claim
that the images of $(b\otimes 1\otimes 1\otimes c)$ under $\tilde\alpha$
constitute a basis of the image of $\tilde\alpha$, i.e. of $D_\alpha$.} They generate
the image, since we have the following equalities:
$\tilde\alpha_1(b\otimes 1\otimes X\otimes c)=\tilde\alpha_1(-b\otimes 1\otimes 1\otimes
Xc)$, $\tilde\alpha_2(b\otimes Y\otimes 1\otimes c)=-\tilde\alpha_2(b\otimes
1\otimes Y\otimes c)$ and $\tilde\alpha_1(b\otimes
Y\otimes X\otimes c)=\tilde\alpha_1(b\otimes Y\otimes 1\otimes Xc)$.
On the other hand we know that
$\mV\theta_s\CC_t\CC_s P(w_0)\cong\mV\theta_sP(w_0)\cong \cC\oplus \cC$, hence the claim follows.

{\it We claim that the images of $(b\otimes 1\otimes
1\otimes c)$ under $\tilde\beta$
constitute a basis of the image of $\tilde\beta$, i.e. of $D_\beta$}. Again,
it is sufficient to show that they generate the image.  Let $\cB$ denote their
$\mC$-span. Note that the $\{d\otimes x\otimes y\otimes c\}$  for $c$ running through
a $\mC$-basis of $\cC$, $x\in\{1,X\}$, and $d$, $y\in\{1,Y\}$ are a $\mC$-basis
of $\cC\otimes_{\cC^t}\cC\otimes_{\cC^s} \cC\otimes_{\cC^t} \cC$. We will
frequently use the following formulas
\begin{eqnarray*}
   \tilde\beta(aY\otimes b\otimes c\otimes d)&=&-\tilde\beta(a\otimes
   bY\otimes c\otimes d)\\
\tilde\beta(a\otimes bX\otimes c\otimes d)&=&-\tilde\beta(a\otimes
   b\otimes cX\otimes d)
\end{eqnarray*}
for any $a$, $b$, $c$, $d\in \cC$ without explicitly referring to them. (The
$i$-th formula follows directly from the corresponding property of $\tilde\beta_i$).

Then the claim follows from the following calculations:
\begin{eqnarray*}
\tilde\beta(X\otimes 1\otimes
1\otimes c)&=&\tilde\beta\big(1\otimes (X+2Y)\otimes 1\otimes c+1\otimes
2Y\otimes 1\otimes c\big)\\
&=&\tilde\beta\big(1\otimes 1\otimes (4Y+8X)\otimes c+1\otimes 1\otimes 7X\otimes
c\big)\\
&=&\tilde\beta\big(1\otimes 1\otimes(15X+30Y)\otimes c-26(1\otimes 1\otimes
Y\otimes c)\big)\\
&=&\tilde\beta\big(1\otimes 1\otimes 1\otimes (15X+30Y)c-26(1\otimes 1\otimes
Y\otimes c)\big)
\end{eqnarray*}
 (we used that $X+2Y$ is $t$-invariant and $Y+2X$ is $s$-invariant). Hence,
 \begin{eqnarray}
\label{11yc}
 \tilde\beta(1\otimes 1\otimes Y\otimes c)\in\cB.
 \end{eqnarray}
This implies that $\tilde\beta(1\otimes X\otimes 1\otimes
c)=-\tilde\beta(1\otimes 1\otimes X\otimes c)\in \cB$ and also
$\tilde\beta(1\otimes X\otimes Y\otimes c)=\tilde\beta(-1\otimes 1\otimes
XY\otimes c)\in\cB$. Therefore, $\tilde\beta(Y\otimes 1\otimes 1\otimes
c)=\tilde\beta(-1\otimes Y\otimes 1\otimes c)=\tilde\beta(-1\otimes 1\otimes
(2X+Y)\otimes c+1\otimes 2X\otimes 1\otimes c)\in\cB$. Finally
$\tilde\beta(Y\otimes X\otimes 1\otimes c)=-1\otimes XY\otimes 1\otimes
c\in\cB$, $\tilde\beta(Y\otimes 1\otimes Y\otimes c)=-\tilde\beta(1\otimes
Y\otimes Y\otimes c)\in\cB$ and $\tilde\beta(Y\otimes X\otimes Y\otimes
c)=-\tilde\beta(1\otimes XY\otimes Y\otimes c)\in\cB$. The claim follows.

Now one can choose a morphism of $\cC$-bimodules
\begin{eqnarray*}
\phi:\cC\otimes_{\cC^t}\cC\otimes_{\cC^s} \cC\otimes_{\cC^t}\cC\rightarrow
\cC\otimes_{\cC^s}\cC\otimes_{\cC^t} \cC\otimes_{\cC^s}\cC
\end{eqnarray*}
which maps $1\otimes 1\otimes 1
\otimes 1$ to  $1\otimes 1\otimes 1
\otimes 1$ and induces an isomorphism on the subbimodules given by inclusions
$i$ and $j$ of $F$. (This choice is possible, since the head
 of $F\Delta(e)$ is simple and isomorphic to $\Delta(w_0)$ ``sitting in
 minimal possible degree'').
In particular, $\phi$ defines a bijection on
the bases constructed above. Hence we constructed an isomorphism $\psi:
\mV\theta_{s}\CC_t\CC_{s}\mV^{-1}\cong\mV\CC_t\CC_{s}\theta_t\mV^{-1}$ giving
rise to an isomorphism
$\theta_{s}\CC_t\CC_{s}\cong\CC_t\CC_{s}\theta_t$ when restricted to the
category of projective-injective objects. By the remarks at the beginning of
the proof we get an
isomorphism of endofunctors on $\cO_0$. This completes the proof.
\end{proof}

\begin{remark}
{\rm Using graded versions of all functors involved (which requires a
further development of some theory  from \cite{MS}, \cite{StGrad}) one could
give a more conceptual proof as follows: One can first show that there is an
embedding of $\theta_s\theta_t\oplus\theta_s$ into $\theta_s\theta_t\theta_s$,
whose cokernel is isomorphic to $\theta_s\CC_t\CC_s$ on the one hand side, but
also to  the quotient $Q$ of the homogeneous inclusion $\theta_s\theta_t\inj
F$ of degree one on the other side. Analogously, there is an embedding
of $\theta_s\theta_t\oplus\theta_t$ into $\theta_t\theta_s\theta_t$, whose
cokernel is isomorphic to $\CC_t\CC_s\theta_t$ on the one hand side, but also
to the functor $Q$. This implies then the first isomorphism of the previous
proposition.
\erem}
\end{remark}

\section*{Acknowledgments}
A part of this research was done during the visit of the
second author to Uppsala  University. The financial support and
hospitality of Uppsala  University are gratefully acknowledged.
The first author was partially supported by STINT, The Royal Swedish
Academy of Sciences and The Swedish Research Council. The second
author was supported by CAALT and EPSRC.

Thanks are due to Henning Haahr Andersen for
several very interesting discussions and comments on previous
versions of the paper. We also would like to thank Wolfgang Soergel
for helpful discussions related to the content of the paper, and to the
referee for valuable comments.

\vspace{0.5cm}

\noindent
Volodymyr Mazorchuk, Department of Mathematics, Uppsala University,
Box 480, 751 06, Uppsala, Sweden, e-mail: {\tt mazor\symbol{64}math.uu.se},
web: {``http://www.math.uu.se/$\tilde{\hspace{1mm}}$mazor/''}.
\vspace{0.2cm}

\noindent
Catharina Stroppel,
Department of Mathematics,
University of Glasgow,
University Gardens,
Glasgow G12 8QW, United Kingdom,
e-mail: {\tt cs\symbol{64}maths.gla.ac.uk}.


\begin{thebibliography}{9999}
\bibitem[An]{Afilt}
{H.H. Andersen}, {\em Filtrations and tilting modules.}
Ann. Sci. \'Ec. Norm. Sup. (4) 30, (1997), 353-366.
\bibitem[AL]{AL}
{H. H. Andersen, N. Lauritzen}, {\em Twisted Verma modules.} in:
Studies in  Memory of Issai Schur, 1-26, v. 210, Progress in Math,
Birkh{\"a}user, Basel, 2002.
\bibitem[AS]{AS}
{H. H. Andersen, C. Stroppel},  {\em Twisting functors on $\mathcal{O}$.}
Represent. Theory 7 (2003), 681-699.
\bibitem[Ara]{Ara}
T.~Arakawa, {\em Vanishing of cohomology associated to quantized Drinfeld-Sokolov
reduction}. Int. Math. Res. Not. 2004, no. 15, 730--767.
\bibitem[Ark]{Ark}
{S.~Arkhipov}, {\em Semi-infinite cohomology of associative algebras and bar
duality}. Internat. Math. Res. Notices 1997, no. 17, 833--863.
\bibitem[AR]{AR} M. Auslander, I. Reiten,
{\em Applications of contravariantly finite subcategories.}
Adv. Math. 86 (1991), no. 1, 111--152.
\bibitem[Bac]{Ba}
E. Backelin, {\em The Hom-spaces between projective functors.}
Represent. Theory 5 (2001), 267--283 (electronic).
\bibitem[Bae]{Baez}
J. Baez, {\em Link invariants of finite type and perturbation theory.}
 Lett. Math. Phys. 26, (1992), (1), 43--51.
\bibitem[BBM]{BBM}
A.~Beilinson, R.~Bezrukavnikov, I.~Mirkovi{\'c}, Tilting exercises.
Mosc. Math. J. 4 (2004), no. 3, 547--557, 782.
\bibitem[Be]{Be}
{I. N. Bernstein}, {\em Trace in categories.} Operator algebras,
unitary representations, enveloping algebras, and invariant
theory (Paris, 1989), 417--423, Progr. Math., 92, Birkh{\"a}user Boston,
Boston, MA, 1990.
\bibitem[BFK]{BFK}
{I. N. Bernstein, I. Frenkel, M. Khovanov},
{\em A categorification of the Temperley-Lieb algebra and
Schur quotients of $U(\mathfrak{sl}_2$) via projective and
Zuckerman functors.} Selecta Math. (N.S.), 5 (1999), no. 2, 199--241.
\bibitem[BG]{BG}
{I. N. Bernstein, S. I. Gelfand}, {\em Tensor products of finite-
and infinite-dimensional representations of semisimple Lie algebras.}
Compositio Math. 41 (1980), no. 2, 245--285.
\bibitem[BGG]{BGG}
{I. N. Bernstein, I. M. Gelfand, S. I. Gelfand}, {\em A certain
category of $\g$-modules.} Funkcional. Anal. i Prilozen. 10 (1976),
no. 2, 1--8.
\bibitem[Bi]{Birman}
J. S. Birman, {\em New points of view in knot theory.},
Bull. Amer. Math. Soc. (N.S.), 28, (1993) no. 2, 253-287.
\bibitem[CR]{CR}
J.~Chuang, R.~Rouquier, {\em Derived equivalences for symmetric groups
and $\mathfrak{sl}_2$-ca\-te\-go\-ri\-fi\-ca\-tion}.  math.RT/0407205,
to appear in Annals Math.
\bibitem[Co]{Corran}
R. Corran, {\em A normal form for a class of monoids including the singular
braid monoids.} J. Algebra, 223, (2000), no.1, 256-282.
\bibitem[DG]{DasbachGemein}
O. Dasbach, B. Gemein, {\em A faithful representation of the singular braid
monoid on three strands.} in: Knots in Hellas '98, Ser. Knots
Everything, no. 24, 48--58
\bibitem[De]{De}
V. Deodhar, {\em On a construction of representations and a problem of
Enright.}  Invent. Math.  57  (1980), no. 2, 101-118.
\bibitem[En]{En}
{T. J. Enright}, {\em On the fundamental series of a real semisimple Lie
algebra: their irreducibility, resolutions and multiplicity formulae.}
Ann. of Math. (2)  110  (1979), no. 1, 1-82.
\bibitem[EW]{EW}
{T. J. Enright, N. R. Wallach},  {\em  On homological algebra and
representations of Lie algebras.} Duke Math. J. 47 (1980), no. 1, 1-15.
\bibitem[GP]{GP}
J. Gonz\'alez-Meneses, L. Paris, {\em Vassiliev invariants for braids on
surfaces.} Trans. Amer. Math. Soc., 356, (2004), no. 1, 219-243.
\bibitem[Ha]{Ha}
D. Happel, {\em Triangulated categories in the representation theory
of finite-dimensional algebras.} London Mathematical Society Lecture
Note Series, 119. Cambridge University Press, Cambridge, 1988.
\bibitem[HR]{HR}
D. Happel, C. M. Ringel, {\em Tilted algebras.}  Trans. Amer. Math.
Soc.  274  (1982), no. 2, 399--443.
\bibitem[Ir1]{I}
{R. Irving}, {\em Shuffled Verma modules and principal series modules
over complex semisimple Lie algebras.} J. London Math. Soc. (2) 48
(1993), no. 2, 263--277.
\bibitem[Ir2]{I2}
{R. Irving}, {\em Projective modules in the category $\Oo\sb S$:
self-duality.} Trans. Amer. Math. Soc.  291  (1985),  no. 2,
701--732.
\bibitem[Ja]{Ja}
J. Jantzen, {\em Einh{\"u}llende Algebren halbeinfacher Lie-Algebren}.
Ergebnisse der Mathematik und ihrer Grenzgebiete (3), 3.
Springer-Verlag, Berlin, 1983.
\bibitem[Jo1]{J}
{A. Joseph}, {\em The Enright functor on the Bernstein - Gelfand -
Gelfand category $\Oo$.} Invent. Math. 67 (1982), no. 3, 423--445.
\bibitem[Jo2]{JCompl}
{A. Joseph}, {\em Completion functors in the ${\mathcal O}$ category.}
Noncommutative harmonic analysis and Lie groups (Marseille, 1982),
80--106, Lecture Notes in Math., 1020, Springer, Berlin, 1983.
\bibitem[KL]{KL}
D. Kazhdan, G. Lusztig, {\em Representations of Coxeter groups and Hecke
algebras}. Invent. Math.  53  (1979), no. 2, 165--184.
\bibitem[Kh]{Kh}
O. Khomenko, {\em Categories with projective functors.} Proc. London Math.
Soc. 90 (2005), no. 3, 711-737.
\bibitem[KM]{KM}
{O. Khomenko, V. Mazorchuk}, {\em On Arkhipov's and Enright's functors}.
Math. Z. 249  (2005),  no. 2, 357--386.
\bibitem[KSX]{KSX}
S. K{\"o}nig, I. Slung{\aa}rd, C.Xi, {\em Double centralizer properties,
dominant dimension, and tilting modules.}  J. Algebra  240  (2001),
no. 1, 393--412.
\bibitem[La]{La}
G. Lallement, {\em Semigroups and combinatorial applications.} Pure and
Applied Mathematics. A Wiley-Interscience Publication. John Wiley \&
Sons, New York-Chichester-Brisbane, 1979.
\bibitem[Ma]{M}
{V. Mazorchuk}, {\em Twisted and shuffled filtrations on tilting modules.}
Math. Reports of the Academy of Science of the Royal Society of Canada
25 (2003), no. 1, 26--32.
\bibitem[MO]{MO}
V. Mazorchuk and S. Ovsienko, {\em Finitistic dimension of properly
stratified algebras.} Adv. Math. 186  (2004),  no. 1, 251--265.
\bibitem[MS1]{MS}
{V. Mazorchuk, C. Stroppel}, {\em
Translation and shuffling of projectively presentable modules and a
categorification of a parabolic Hecke module}. Trans. Amer. Math. Soc.
357  (2005), 2939-2973.
\bibitem[MS2]{KMS}
{V.~Mazorchuk, C.~Stroppel}, {\em Projective-injective modules,
Serre functors and symmetric algebras}, math.RT/0508119,
to appear in J. Reine Angew. Math.
\bibitem[Or]{Orekov}
S. Orevkov: {\em Solution of the word problem in the singular braid group.}
Turkish J. Math., 28, (2004), no. 1, 95--100.
\bibitem[Rin]{R}
{C. M. Ringel}, {\em The category of modules with good filtrations over
a quasi-hereditary algebra has almost split sequences.} Math. Z. 208
(1991), no. 2, 209--223.
\bibitem[Ric]{Rickard}
J. Rickard, {\em Translation functors and equivalences of derived
categories for blocks of algebraic groups.} Finite-dimensional
algebras and related topics (Ottawa, ON, 1992), 255--264, NATO
Adv. Sci. Inst. Ser. C Math. Phys. Sci., 424, Kluwer Acad. Publ.,
Dordrecht, 1994.
\bibitem[Re]{Reiten}
I. Reiten, {\em Tilting theory and homologically finite subcategories
with applications to quasihereditary algebras.} to appear in the
Tilting Handbook edited by L. Angeleri H\"ugel, D. Happel and H. Krause.
\bibitem[So1]{SoInv}
W. Soergel, {\em {$\mathfrak n$}-cohomology of simple highest weight modules on
              walls and purity}, Invent. Math., 98, (1989), no. 3, 565--580.
\bibitem[So2]{S}
{W. Soergel}, {\em Kategorie $\Oo$, perverse Garben und Moduln
{\"u}ber den Koinvarianten zur Weylgruppe.} J. Amer. Math. Soc.
3 (1990), no. 2, 421--445.
\bibitem[So3]{SHC}
{W. Soergel}, {\em The combinatorics of Harish-Chandra bimodules.}
J. Reine Angew. Math.  429  (1992), 49--74.
\bibitem[So4]{So4}
{W. Soergel}, {\em Character formulas for tilting modules over Kac-Moody algebras}.
Represent. Theory  2  (1998), 432--448.
\bibitem[St1]{StGrad}
C. Stroppel, {\em Category $\Oo$: Gradings and Translation Functors.}
J. Algebra 268 (2003), no. 1, 301-326.
\bibitem[St2]{St4}
C. Stroppel, {\em Homomorphisms and extensions of principal series
representations.}  J. Lie Theory  13  (2003),  no. 1, 193--212.
\bibitem[St3]{Stquiv}
C. Stroppel, {\em Category $\cO$: Quivers and Endomorphism rings of
  Projectives.} Represent. Theory 7 (2003), 322--345.
\bibitem[Ve]{Ve}
 V. Vershinin, {\it Vassiliev invariants and singular braids},
 Russ. Math. Surv., 1998, 53 (2), 410--412.
\bibitem[Wa]{Wa}
T. Wakamatsu, {\em Stable equivalence for self-injective algebras and
a generalization of tilting modules}.  J. Algebra  134  (1990),
no. 2, 298--325.
\bibitem[Zu]{Z}
G. Zuckerman, {\em Tensor products of finite and infinite dimensional
representations of semisimple Lie groups}, Ann. Math. (2) 106 (1977), no. 2, 295--308.
\end{thebibliography}
\end{document}